\documentclass{amsart}

\usepackage{amssymb}
\usepackage[all]{xy}
\usepackage{hyperref}

\newtheorem{thm}{Theorem}

\newtheorem{lem}[thm]{Lemma}
\newtheorem{cor}[thm]{Corollary}

\newtheorem{prop}[thm]{Proposition}

\theoremstyle{definition}
\newtheorem{defn}[thm]{Definition}

\newtheorem{say}[thm]{}
\newtheorem{exmp}[thm]{Example}

\newtheorem{ques}[thm]{Question}    

\newtheorem{rem}[thm]{Remark}          

\newtheorem*{ack}{Acknowledgments}      

\newtheorem{defn-thm}[thm]{Definition--Theorem}  
\newtheorem{defn-lem}[thm]{Definition--Lemma}  

\theoremstyle{remark}

\renewcommand{\c}[0]{{\mathbb C}}  

\renewcommand{\o}[0]{{\mathcal O}} 
\newcommand{\z}[0]{{\mathbb Z}}
\newcommand{\n}[0]{{\mathbb N}}

\renewcommand{\a}[0]{{\mathbb A}}

\newcommand{\p}[0]{{\mathbb P}}
\newcommand{\f}[0]{{\mathbb F}}
\newcommand{\q}[0]{{\mathbb Q}}
\newcommand{\map}[0]{\dasharrow}
\newcommand{\qtq}[1]{\quad\mbox{#1}\quad}
\newcommand{\spec}[0]{\operatorname{Spec}}

\newcommand{\mult}[0]{\operatorname{mult}}

\newcommand{\supp}[0]{\operatorname{Supp}}

\newcommand{\im}[0]{\operatorname{im}}

\newcommand{\aut}[0]{\operatorname{Aut}}

\newcommand{\sing}[0]{\operatorname{Sing}}

\newcommand{\jac}[0]{\operatorname{Jac}} 

\newcommand{\chr}[0]{\operatorname{char}}

\newcommand{\rdown}[1]{\lfloor{#1}\rfloor}

\newcommand{\ord}[0]{\operatorname{ord}}

\newcommand{\tsum}[0]{\textstyle{\sum}}

\def\into{\DOTSB\lhook\joinrel\to}

\def\loccoh#1.#2.#3.#4.{H^{#1}_{#2}(#3,#4)}

\DeclareMathAlphabet{\mathchanc}{OT1}{pzc}%
                                {m}{it}

\newcommand{\gm}[0]{{\mathbb G}_m}

\usepackage[all]{xy}\xyoption{dvips}

\newcommand{\norm}[0]{\operatorname{norm}}

\newcommand{\find}[0]{\operatorname{f-ind}}
\newcommand{\fdiv}[0]{\operatorname{div}}
\newcommand{\mfrd}[0]{\operatorname{\mathfrak{d}}}

\begin{document}
\bibliographystyle{amsalpha}
\hfill\today

 \title{Pell surfaces}
 \author{J\'anos Koll\'ar}

\begin{abstract}  In 1826 Abel started the study of the polynomial Pell equation
$x^2-g(u)y^2=1$. Its solvability in polynomials $x(u), y(u)$ depends on a certain torsion point on the Jacobian of the hyperelliptic curve $v^2=g(u)$.
In this paper we study the affine surfaces  defined by the Pell equations in 3-space with coordinates $x, y,u$, and aim to describe all affine lines on it. These are polynomial solutions of the equation  $x(t)^2-g(u(t))y(t)^2=1$. Our results are rather complete when the degree of $g$ is even but the odd degree cases are left completely open. For even degrees we also describe all curves on these Pell surfaces that have only 1 place at infinity. 
\end{abstract}

 \maketitle

The classical Pell equation is  $x^2-dy^2=1$; its rational solutions correspond to the units in the number field  $\q\bigl(\sqrt{d}\bigr)$.
 Abel  \cite{abel1826} studied the {\it polynomial  Pell equation}  $x^2-g(u)y^2=1$, where  $g(u)$ is a polynomial, looking for solutions where $x=x(u), y=y(u)$ are also polynomials in $u$. 
In this note we look at the polynomial  Pell equation as an affine algebraic surface over a field $k$
$$
S_g:=\bigl(x^2-g(u)y^2=1\bigr)\subset \a^3_{xyu},
\eqno{(*)}
$$
and aim to describe all non-constant morphisms
$\a^1\to S_g$. Equivalently, all solutions of the equation
$$
x(t)^2-g\bigl(u(t)\bigr)y(t)^2=1\qtq{where} x(t), y(t), u(t)\in k[t].
\eqno{(**)}
$$
Solutions of the polynomial  Pell equation correspond to those
$\bigl(x(t), y(t), u(t)\bigr)$ 
for which $u(t)=t$; these are the sections of the coordinate projection
$\pi:S_g\to \a^1_u$.

If $g(u)=u$ then solving ($**$) is equivalent to solving all
polynomial Pell equations simultaneousy, thus we focus on the cases when  $\deg g\geq 2$.  

\begin{defn}[Affine lines in  varieties] \label{aff.lin.aff.var.say}
 Let $X$ be a quasi-affine variety. We call a closed curve  $B\subset X$ an
{\it affine line} if $B\cong \a^1$ and a 
{\it singular affine line} if the normalization of $B$ is isomorphic to  $\a^1$

Let $\phi:\a^1\to X$ be a non-constant morphism and  $\phi_C:C\to X$  the normalization of
the closure of its image. Then $\phi$ lifts to $\tau_C:\a^1\to C$. 
Thus  $C\cong \a^1$, hence $\phi(\a^1)$ is a singular affine line.
  Thus understanding all   non-constant morphisms  $\phi:\a^1\to X$\
is equivalent to  understanding all   (possibly singular) affine lines  $B\subset X$. 

\medskip
{\it Obvious affine lines \ref{aff.lin.aff.var.say}.1.} Every Pell surface  $S_g$ contains $\leq 2\deg g+2$  obvious affine lines.
For every root $g(c)=0$ we have 2 {\it vertical}  affine lines $t\mapsto (x=\pm 1, y=t, u=c)$ and we also have the 2 {\it trivial sections}  $t\mapsto (x=\pm 1, y=0, u=t)$.  
\end{defn}

Our first observation is  that, in many cases, the theory of  polynomial  Pell equations gives all affine lines on  Pell surfaces.

\begin{thm}\label{aff.lin.sec.thm} Let $k$ be a perfect field of characteristic $\neq 2$ 
 and  $g(u)\in k[u]$ a  polynomial of even degree. Then every 
(possibly singular) affine line on  the Pell surface
$$
S_g:=\bigl(x^2-g(u)y^2=1\bigr)\subset \a^3_{xyu},
$$
is either vertical or a section of the coordinate projection
$\pi:S_g\to \a^1_u$. 
\end{thm}

By contrast, if $\deg g$ is odd then there are no nontrivial sections, but there can be other singular affine lines on $S_g$. 
A discussion of the degree 3 case, due to Zannier, is given in Example~\ref{deg.3.exmp}

\medskip

In the theory of projective surfaces, lot of attention has been paid to
understanding rational curves on surfaces of Kodaira dimension 0 or 1. 
 Pell surfaces are affine analogs of elliptic K3 surfaces  (the $\deg g=2$ case) and of properly elliptic surfaces    $S\to \p^1$ (the $\deg g\geq 3$ cases). 
For elliptic surfaces the Mordell-Weil group describes the sections; these are quite well understood, see \cite{ss-mw}.  Elliptic K3 surfaces  usually contain  infinitely many other  rational curves; these are not well understood. 
For properly elliptic surfaces  one expects very few other rational curves, but I do not know a single example of a projective surface defined over $\bar \q$,  with Kodaira dimension $\geq 0$  and containing infinitely many rational curves that are all explicitly known. See  \cite{MR3606999} or 
Example~\ref{only.sects.from K3} for some examples over $\c$.

The analogous question on open surfaces is the description of 
affine lines on  surfaces of log Kodaira dimension 0 or 1. By Proposition~\ref{pell.surf.lem.1} the log Kodaira dimension of $S_g$ is 0 if $\deg g=2$ and 1  if $\deg g\geq 3$. Thus, combining Theorem~\ref{aff.lin.sec.thm}
with known results on solutions of the  polynomial  Pell equation, 
we get examples of open surfaces with log Kodaira dimension $\geq 0$  containing infinitely many  affine lines, all of which one can list
 explicitly. 
We discuss in detail the
simplest case $S_2:=\bigl(x^2-(u^2-1)y^2=1\bigr) $ in Example~\ref{u2-1.exmp}.

A complete list of all  surfaces of log Kodaira dimension 0 that contain infinitely many affine lines is given in \cite{MR3649236}. 
However, an enumeration  of all   affine lines seems to be known only
for the Pell surface in Example~\ref{u2-1.exmp}.
A much studied example is  $T:=\p^2\setminus(\mbox{smooth cubic})$.
Affine lines of degree $\leq 7$ in $T$ are determined in
\cite{takahashi}. Almost all enumerative invariants involving rational curves in $T$ are computed in 
\cite{MR2435425}, but the number of affine lines of a given degree is left undetermined.

\begin{defn} \label{1.pl.at.inf.defn}
Let $D$ be an affine curve over a field $k$ with normalization $D^n$ and
 smooth compactification  $D^n\subset \bar D$.
The geometric points of $\bar D\setminus D^n$ are the {\it places at infinity} 
of $D$.  For example, over $\c$ the curve   $\bigl(y^n=g(x)\bigr)$ has only 1 place at infinity iff
$(n, \deg g)=1$. 
If $D$ has only 1 place at infinity then the invertible regular functions on $D$ are constants. 

Given an open surface $S$, it is of interest to study curves
$D\subset S$ that have only 1 place at infinity. There are especially complete results about  $S=\a^2$; see  \cite{MR0338423, MR0379502, MR1697368, MR1951527, MR2477978, MR3089043, de-du} and the references there.

Surfaces with log Kodaira dimension 0 sometimes contain 
no curves with  only 1 place at infinity
(for example  $\c^*\times \c^*$); in other cases they contain 
positive dimensional families of such curves. For example,
let $E\subset \p^2$ be a smooth cubic and $L$ a flextangent.
Then every member of the linear system  $|E, 3L|$ meets $E$ only at the flex.
Thus we get a 1-dimensional family of curves in $\p^2\setminus E$ that have
geometric genus 1 and only 1 place at infinity. 
We get  larger genus examples starting with higher order torsion points on $E$.
\end{defn}

Pell surfaces also give examples that contain infinitely many affine lines but no other curves with only 1 place at infinity.

\begin{thm}\label{1.pl.sec.thm} Let $k$ be a field of characteristic $\neq 2$ 
 and  $g(u)\in k[u]$ a nonzero polynomial of even degree. Then every 
curve  with only 1 place at infinity  on  the Pell surface
$S_g:=\bigl(x^2-g(u)y^2=1\bigr)\subset \a^3_{xyu}$
is an affine line.
\end{thm}

By contrast, every odd degree Pell surface contains infinitely many curves with only 1 place at infinity, see Example~\ref{1.pl.odd.deg.exmp}.

\medskip

In writing this article, I tried to build the technical machinery only as needed and start with  elementary treatments of significant special cases whenever possible. 

 Section~\ref{sec.1} discusses thes simplest examples and their application to the undecidability of the embedding problem for affine varieties, due to
Chilikov and Kanel-Belov \cite{kb-ch}, which needs only the knowledge of affine lines on
the  simplest Pell  surface  
$S_2:=\bigl(x^2-(t^2-1)y^2=1\bigr)$. 

 Section~\ref{sec.2} gives the proof of Theorem~\ref{aff.lin.sec.thm} using Abel's method and the rest of the paper is devoted to proving Theorem~\ref{1.pl.sec.thm}.

A general introduction to the geometry of Pell  surfaces is in  Section~\ref{sec.3}, followed by an introduction to polynomial Pell equations in  Section~\ref{sec.4}. This is mostly based on   \cite{scherr} and \cite{zan1, zan2}. We quickly revisit Hazama's approach to  Theorem~\ref{aff.lin.sec.thm} in Section~\ref{sec.5}. 

In Section~\ref{sec.6} we reduce Theorem~\ref{1.pl.sec.thm} to a question about  maps between the first  homology groups of certain non-compact algebraic curves. This turns out to be easy using topology, which leads to a  proof of Theorem~\ref{1.pl.sec.thm}  over $\c$.
This is discussed in    Sections~\ref{sec.7}--\ref{sec.8}.

The positive characteristic case is more complicated.
The  plan of the proof is outlined in Section~\ref{sec.9}, with the details in 
Sections~\ref{pi1.surj.sec}--\ref{fund.ind.def.sec}.

In Section~\ref{endom.sec} we describe all endomorphisms of Pell surfaces.
The hardest case is $S_2=\bigl(x^2-(u^2-1)y^2=1\bigr)$, where our computations rely on the complete enumeration of all affine lines; see  Paragraphs~\ref{end.S2.char.0.say}--\ref{end.S2.char.p.say} for details.

\begin{ack} 
I thank A.A.~Chilikov and A.J.~Kanel-Belov for posing the original question,  D.~Gabai and 
Z.~Scherr for help with  the literature, L.~Chen, S.~Kov\'acs,  M.~Lieblich, B.~Totaro  and J.~Waldron for helpful conversations  and Umberto~Zannier
for many comments, corrections and examples. 
Partial  financial support    was provided  by  the NSF under grant numbers
 DMS-1362960 and  DMS-1440140 while the author was in residence at
MSRI during the Spring 2019 semester.
\end{ack}

\section{Examples and applications}\label{sec.1}

\begin{exmp} \label{u2-1.exmp}
Let $k$ be a field of characteristic $\neq 2$.
On the Pell surface  
$$
S_2:=\bigl(x^2-(t^2-1)y^2=1\bigr)\subset \a^3_{xyt}
\eqno{(\ref{u2-1.exmp}.1)}
$$
every (possibly singular) affine line  is smooth. 
Besides the $\leq 6$ obvious ones listed in (\ref{aff.lin.aff.var.say}.1), we
immediately see the solution  $x=t, y=1$. As with the usual Pell equation, we then get other solutions by the formula
$$
x_n(t)+y_n(t)\sqrt{t^2-1}=\bigl(t+\sqrt{t^2-1}\bigr)^n.
$$ 
Thus we  have the  infinite sequence of sections $\Sigma_n$ given by
$$
\begin{array}{lcl}
x_n(t)& = & \sum_{i=0}^{\rdown{n/2}} \binom{n}{2i}t^{n-2i}(t^2-1)^i,\\
y_n(t)& = & \sum_{i=0}^{\rdown{n/2}} \binom{n}{2i+1}t^{n-2i-1}(t^2-1)^i,\\
u_n(t)& = & t,
\end{array}
\eqno{(\ref{u2-1.exmp}.2)}
$$
for $n\geq 1$, and also   $ \bigl(\pm x_n(t),  \pm y_n(t), t\bigr)$ for all sign choices.  \cite{MR0491583} proves that these give all solutions, though this was most likely already known to Abel. The  paper \cite[4.3]{hazama}  shows that there are no affine lines
$\bigl(x(s), y(s), t(s)\bigr)$ for which $t^2(s)-1$ has only simple roots.
We see in Section~3  that the latter restriction is not necessary.

The intersection points of the affine lines on $S_2$ have remarkable properties.

(\ref{u2-1.exmp}.3) For $t=1$ only the $i=0$ summands in (\ref{u2-1.exmp}.2) are nonzero, thus 
we obtain---as observed by \cite{MR0491583}---that the affine lines on $S_2$   intersect
the line $L=(x-1=u-1=0)$ precisely at the points
$(1,n,1)$ where $n\in \z\setminus\{0\}$. 
In \cite{MR0491583} this was used  to prove 
that there are some  undecidable questions in algebraic geometry.

(\ref{u2-1.exmp}.4) Observe that $x_n(t)$ and $y_n(t)$ are the Chebyshev polynomials of the first and second kind, defined by the properties
$$
T_n(\cos\theta)=\cos(n\theta)\qtq{and} U_{n-1}(\cos\theta)={\sin(n\theta)}/{\sin\theta}.
$$
The identity  
$$
\cos^2(n\theta)-(\cos^2\theta-1){\sin^2(n\theta)}/{\sin^2\theta}=
\cos^2(n\theta)+\sin^2(n\theta)=1
$$
shows that they lie on the surface $S_2$.

We see in Claim~\ref{int.of.sec.S2}.2 that  the projection of all intersection points of all sections to $\a^1$ is the set
$R_{\infty}:=\{\cos(2\pi \alpha): \alpha\in \q\}$.

(\ref{u2-1.exmp}.5) Let $\phi:S_2\to S_2$ be a dominant endomorphism.
Then $\phi$ has only finitely many exceptional curves, 
hence all but finitely many of the affine lines on $S_2$ are mapped to
affine lines on $S_2$. Since we have a good description of all affine lines on $S_2$,
we can use them to  determine all automorphisms and endomorphisms of $S_2$.
For the other Pell surfaces this turns out to be much easier;
 see Theorem~\ref{endom.thm} for the precise statement.
\end{exmp}

The next  application of Example~\ref{u2-1.exmp} gave the original motivation to consider this question. Its  proof uses both the explicit description of the sections  (\ref{u2-1.exmp}.2) and the fact that there are no other affine lines on $S_2$. 
 (We call a map $\phi:Y\to X$ {\it non-degenerate} if  $\dim Y=\dim\bigl(\phi(Y)\bigr)$.) 

\begin{thm} \cite[Thm.4]{kb-ch}  
For affine varieties  $X$ defined over $\q$,
 the following questions are all  algorithmically undecidable.
\begin{enumerate}
\item Is there a closed embedding $\a^{11}\to X$ defined over $\q$?
\item Is there a non-degenerate morphism $\a^{11}\to X$ defined over $\q$?
\item Is there a closed embedding $\a^{11}\to X_{\c}$ defined over $\c$?
\item Is there a non-degenerate morphism $\a^{11}\to X_{\c}$ defined over $\c$?
\end{enumerate}
\end{thm}

\begin{rem} The use of $\a^{11}$ is almost certainly an artifice of the proof and there are probably many other algorithmically undecidable questions in algebraic geometry. Roughly speaking, undecidability could occur every time a property holds for certain objects that correspond to points in a  countably infinite union of  
subvarieties in a moduli space.
\end{rem}

\begin{exmp} \label{u2-c.exmp}
Let $k$ be a field of characteristic $\neq 2$. The general degree 2 Pell surface can be written as
$$
S_{ac}:=\bigl(x^2-(at^2-c)y^2=1\bigr)\subset \a^3_{xyt},
\eqno{(\ref{u2-c.exmp}.1)}
$$
where $a, c\neq 0$. We note in Claim~\ref{review.ppe}.2 that there are 
no solutions if $a$ is not a square; thus we may as well assume that $a=1$.
Following the solution of (\ref{u2-1.exmp}.1), we get
$$
x=t\tfrac{1}{\sqrt{c}},\quad y=\tfrac1{\sqrt{c}}
$$
as a solution. This is in $k[t]$ iff $ c\in k^2$.
We take  its square 
$$
\bigl(t\tfrac{1}{\sqrt{c}}+\sqrt{c}\sqrt{t^2-c}\bigr)^2=
\bigl(t^2\tfrac{1}{c}+c(t^2-c)\bigr)+ 2t\sqrt{t^2-c}.
$$
to get the $k[t]$-solution
$$
x_2=\bigl(\tfrac1{c}+c\bigr)t^2-c^2, \quad y_2=2t.
$$
The  other $k[t]$-solutions are given by the formula
$$
\pm x_{2n}(t)\pm y_{2n}(t)\sqrt{t^2-c}=\bigl(x_2+y_2\sqrt{t^2-c}\bigr)^n.
$$ 
\end{exmp}

The following example of cubic Pell surfaces was explained to me by Zannier.

\begin{exmp}\label{deg.3.exmp}
 Let $S_g$ be a Pell surface
$(x^2-g(u)y^2=1)$ where $g(u)$ is a cubic with simple roots.
Then $S_g$ has no sections but we claim that it has 
 infinitely many infinite families of double sections.

To see this note that,
as we discuss in Paragraph~\ref{high.tors.exmp}, there are infinitely many different constants
$c$ for which  $x^2-(u{-}c)g(u)y^2=1$  has  nontrivial solutions.
If $\bigl(x_c(u),y_c(u)\bigr)$ is such a  solution then
$$
x_c(t^2+c)^2-g(t^2+c) \bigl(t y(t^2+c)\bigr)^2=1
$$
shows that    $t\mapsto  \bigl(x_c(t^2+c), ty_c(t^2+c), t^2+c\bigr)$
is a double section. Each value of $c$ yields infinitely many
double sections for which the projection to $\a^1_u$ ramifies over $u=c$. 
Thus different values of $c$ give different double sections. 

See Paragraph~\ref{deg.3.exmp.exp} for more details.
\end{exmp}

\begin{exmp}\label{all.mult.roots}
Let $k$ be a perfect field of characteristic $p\neq 2$. 
If $x^2-g(u)y^2=1$ has a nontrivial solution in $k[u]$ then either $x(u)$ is a $p$th power or $g(u)$ must have at least 2 simple roots. This follows  from the Mason-Stothers theorem. A direct argument is the following.

We claim that  $x(u)^2-1$ has  at least 2 simple roots if
the derivative $x'(u)$ is not identically 0.
To see this write 
 $x(u)^2-1=c\prod (u-a_i)^{m_i}$ and set $h(u):=\prod (u-a_i)^{m_i-1}$. Note that
$h(u)$ divides the derivative if $x(u)^2-1$, which is $2x(u)x'(u)$,  but
it is relatively prime to $x(u)$. So $h(u)$ divides  $x'(u)$. Hence
$$
\tsum_i (m_i-1)=\deg h(u)\leq \deg x'(u)\leq  -1+\tfrac12 \tsum_i m_i.
$$ 
This rearranges to 
$\tsum_i (m_i-2)\leq -2$.

If $x(u)^2-1$ has  exactly 2 simple roots, then, after a linear change of variables we may assume these to be $\pm 1$. Thus
we  have  $x(u)^2-1=(u^2-1)y(u)^2$ for some polynomial $y(u)$. That is, the pair
$\bigl(x(u), y(u)\bigr)$ is one of the solutions of
the Pell equation discussed in  Example~\ref{u2-1.exmp}. 

\end{exmp}

The following is an example of a smooth, projective, elliptic surface over $\p^1$
that contains infinitely many sections but no other rational curves.
See \cite{MR3606999} for a similar result in case of 1 section. 

\begin{exmp}\label{only.sects.from K3}
Let $\pi:X\to \p^1_{st}$ and  $q:\p^1_{uv}\to \p^1_{st}$ be  morphisms.
By base change we get  $$
\pi_q: X_q:=X\times_{\p^1_{st}}\p^1_{uv}\to \p^1_{uv}.
$$
Every section $\sigma:\p^1_{st}\to X$ gives a section of
$\pi_q$, but a rational multi-section usually gives a non-rational multisection.

To understand this, let $C\cong \p^1$ and  $r:C\to \p^1_{st}$.
The fiber product  $C\times_{\p^1_{st}}\p^1_{uv}$ is a curve of bidegree
$(\deg r, \deg q)$ on $C\times\p^1_{uv}$, hence its arithmetic genus is
$(\deg r -1)( \deg q-1)$. If the branch loci of $r$ and $q$ are disjoint, then
$C\times_{\p^1_{st}}\p^1_{uv}$ is smooth, hence non-rational if
$\deg r, \deg q\geq 2$.

To see  more concrete examples, 
let $C_i=(g_i=0)\subset \p^2$ be 2 plane cubics. Assume that they intersect in 9 distintct points $P_1,\dots, P_9$. These are the base points of the pencil of cubics  $\lambda g_1+\mu g_2=0$. By blowing them up  we get a rational elliptic surface
$$
S=(sg_1+tg_2=0)\subset \p^2_{xyz}\times \p^1_{st},
$$
with the 9 exceptional curves giving 9 sections. 
The group of all sections  (usually called the Mordell-Weil group)  is isomorphic to $\z^8$ iff every member of this pencil is irreducible, which holds if no 3 of the 9 intersection points  $(g_1=g_2=0)$ are on a line.

By a double cover of the base we get an elliptic K3 surface
$$
S_2:=\bigl((s^2+t^2)g_1(x,y,z)+(s^2-t^2)g_2(x,y,z)=0\bigr)\subset \p^2_{xyz}\times \p^1_{st}.
$$
This has only countably many rational curves, all defined over $\bar \q$.
In particular, for all of them the branch points of the projection to
$\p^1_{st}$ have algebraic coordinates.  Let now $c$ be any
transcendental number. The branch locus of
$(u,v)\mapsto \bigl(c(u^2+v^2), u^2-v^2\bigr)$ is  $\{(c{:}1), ({-}c{:}1)\}\subset \p^1_{st}$. 
Setting $s=c(u^2+v^2), t=u^2-v^2$ we get a properly elliptic surface
$$
S_3:=\bigl(h_1(u,v)g_1(x,y,z)+ h_2(u,v)g_2(x,y,z)=0\bigr)\subset \p^2_{xyz}\times \p^1_{uv}
$$
where $h_1=c^2(u^2+v^2)^2+(u^2-v^2)^2$ and 
$h_2=c^2(u^2+v^2)^2-(u^2-v^2)^2$.
The group of sections of $S_3\to \p^1_{uv}$ is isomorphic to $\z^8$ 
and, as we noted above, every  rational curve on $S_3$ is either vertical or a section.  

The following claim, whose proof is left as an exercise, allows one to get many concrete examples.

\medskip

{\it Claim \ref{only.sects.from K3}.1.} Pick $c_1,\dots, c_9\in k$ and let
$p_i=(c_i, c_i^3)$ be 9 points on the cubic  $(g_1=0)$ where $g_1:=y-x^3$. 
\begin{enumerate}
\item[(a)] These   9 points are cut out by another cubic $(g_2=0)$ iff
$c_1+\cdots + c_9= 0$.
\item[(b)] Every member of the pencil   $(sg_1+tg_2=0)$ is irreducible
iff no 3 of the $c_i$ sum to 0. 
\end{enumerate}
\end{exmp}

\begin{exmp}\label{1.pl.odd.deg.exmp}
 For any polynomial $h(u)$, the intersection of the Pell surface
$S_g:=\bigl(x^2-g(u)y^2=1\bigr)$ with the surface
$y=h(u)$ is the curve
$$B_{gh}:=\bigl(x^2=g(u)h(u)^2+1\bigr).$$ This curve is hyperelliptic and has only 1 place at infinity iff $\deg g$ is odd.  Thus Theorem~\ref{1.pl.sec.thm} is sharp.
\end{exmp}

\begin{say}[Bogomolov's question] Let $K$ denote either $\bar \f_p$ or $\bar\q$. 
Bogomolov suggested in 1981 that every $K$-point of a K3 surface $S$ might be contained in a rational curve lying on $S$. For Kummer surfaces over $\bar \f_p$ this was proved in \cite{bog-tsc1, bog-tsc2}. 

As a natural analog, one might ask if every  $K$-point of a Pell surface $S$ might be contained in an affine line lying on $S$. This is clearly not the case; we never get all $K$-points on any fiber.

However, if we fix a  Pell surface $S$ over $\f_q$, a quick computation suggests that  affine lines cover a positive proportion of the
$\f_{q^n}$ points of $S$ for every $n$. It would be interesting to understand this better.

\end{say}

\section{Abel's method of continued fractions}\label{sec.2}

\begin{say} Let $k((u^{-1}))$ denote the Laurent series field in $u^{-1}$.
Its elements are of the form $\phi(u)=\tsum_{i\leq N}c_iu^i$ for some $N\in \z$.  
For any  $\phi=\tsum_{i\leq N}c_iu^i\in k((u^{-1}))$  define the {\it polynomial} or {\it integral} 
part of $\phi$ as 
$$
\rdown{\phi}:=\tsum_{i\geq 0}c_iu^i.
$$
Following Abel, the {\it continued fraction expansion} of $\phi$ is defined as follows. Set    $\phi_0:=\phi$. If $\phi_i$ is already defined then
we set
$$
a_i:=\rdown{\phi_i}\qtq{and} \phi_{i+1}:=(\phi_i-a_i)^{-1}=\bigl(\phi_i-\rdown{\phi_i}\bigr)^{-1}.
$$
This represents  $\phi$ as  an infinite {\it continued fraction}
$$
\phi=a_0+\cfrac{1}{a_1+\cfrac{1}{a_2+\cdots}}.
$$
For finite or infinite continued fractions we use the compressed notation 
$$
[a_0,\dots, a_n]:=a_0+\cfrac{1}{a_1+\cfrac{1}{a_2+\cdots}}.
$$
For an infinite continued fraction $\phi=\bigl[a_0(u), a_1(u),\dots\bigr]$,
the $\bigl[a_0(u), \dots, a_n(u)\bigr]$ are called its 
{\it convergents.}  In a precise sense, the convergents give the best approximation of  $\phi$ by rational function; this is called {\it Pad\'e approximation;}   see \cite{pade-wiki}.
(Frequently one writes $\bigl[a_0(u), \dots, a_n(u)\bigr]=p_n(u)/q_n(u)$
where $p_n(u), q_n(u)$ are relatively prime and the pair 
$\bigl(p_n(u), q_n(u)\bigr)$ is called the $n$th  convergent.)

If $g(u)=c_mu^m+\cdots + c_0$ has even degree and $c_m$ is a square in $k$ then we have  a Laurent series expansion in $k((u^{-1}))$ 
$$
\sqrt{g(u)}=\sqrt{c_m}u^{m/2}\sqrt{1+(c_{m-1}/c_m)u^{-1}+\cdots +(c_0/c_m)u^{-m}}.
$$
Thus we get a continued fraction expansion
$$
\sqrt{g(u)}=:\bigl[a_0(u), a_1(u),\dots\bigr].
$$
\end{say}

The following is essentially due to \cite{abel1826}; see also \cite{cheb}
and \cite[Lem.6]{schm} for a complete modern proof.

\begin{thm}\label{abel.cont.frac}
 Let $g(u)$ be a polynomial of even degree and write
$$\sqrt{g(u)}=\bigl[a_0(u),a_1(u), \dots ]$$ as an infinite continued fraction.
Then for every solution  $x_i(u), y_i(u)$ of the Pell equation
$x^2-g(u)y^2=1$, the quotient  $\pm x_i(u)/y_i(u)$ is among the convergents  
$[a_0(u),a_1(u), \dots , a_n(u)]$ for a suitable choice of the sign $\pm$. \qed
\end{thm}

\begin{say}[Proof of Theorem~\ref{aff.lin.sec.thm}]\label{ab-ch.setup.say}
Let $\phi(u)\in k((u^{-1}))$  be a Laurent series with continued fraction expansion
$$
\phi(u)=\bigl[a_0(u), a_1(u),\dots\bigr].
$$
We claim that the continued fraction expansion of $\phi(q(t))$ is given by
$$
\phi(q(t))=\bigl[a_0(q(t)), a_1(q(t)),\dots\bigr].
$$
Due to the inductive definition of the $a_i(u)$, it is enough to show that
$$
\rdown{\phi(q(t))}=\rdown{\phi(u)}\circ q(t).
$$
This needs to be checked for each $u^i$. If $i\geq 0$ then
clearly  $\rdown{q(t)^i}=q(t)^i$ and if $j>0$ then
$$
q(t)^{-j}=c_m^{-j}t^{-mj}\frac{1}{1+(c_{m-1}/c_m)t^{-1}+\cdots +(c_0/c_m)t^{-m}}
$$
shows that $\rdown{q(t)^{-j}}=0$.
Thus every convergent of $\phi(q(t))$ is of the form
$$
\bigl[a_0(u), a_1(u),\dots, a_n(u)\bigr]\circ q(t).
$$
By Theorem~\ref{abel.cont.frac}, every solution 
$\bigl(X(t), Y(t)\bigr)$ of  $x^2-g(q(t))y^2=1$ can be written
as 
$$
\tfrac{X(t)}{Y(t)}=\pm \bigl[a_0(q(t)), a_1(q(t)),\dots, a_n(q(t))\bigr]
$$
for some $n$,  up to sign. Now write
$$
\bigl[a_0(u), a_1(u),\dots, a_n(u)\bigr]=\tfrac{x(u)}{y(u)},
$$
where $x(u),y(u)$ are reatively prime. Then 
$x(q(t)), y(q(t))$ are also reatively prime, so
$X(t)=x(q(t)),\ Y(t)=y(q(t))$, up to multiplicative constants.
We are done since   $\bigl(x(u), y(u)\bigr)$ is a solution of
$x^2-g(u)y^2=1$ iff  $\bigl(x(q(t)), y(q(t))\bigr)$ is a solution of
$x^2-g(q(t))y^2=1$. \qed
\end{say}

\section{Geometry of Pell surfaces}\label{sec.3}

\begin{defn}[Affine Pell surfaces]\label{pell.surf.defn} 
Let $k$ be a field of characteristic $\neq 2$,
$\bar B$ a smooth projective curve over $k$ and  $b_{\infty}\in \bar B(k)$
a point. Set $B:=\bar B\setminus\{b_{\infty}\}$. 
For $g\in k[B]$ let $\deg g$ denote the order of its pole at $b_{\infty}$. 

Let  $g\in k[B]$  be  non-constant. 
We define the corresponding affine {\it Pell surface} as
$$
S_g:=\bigl(x^2-gy^2=1\bigr)\subset \a^2_{xy}\times B.
\eqno{(\ref{pell.surf.defn}.1)}
$$
$S_g$ is smooth and so is the projection $\pi:S_g\to B$.
Moreover,  $\pi:S_g\to B$ is a group scheme with identity section 
$E=(1,0)$ and multiplication
$$
(x_1, y_1)\cdot (x_2, y_2)\mapsto  \bigl(x_1y_1+gx_2y_2, x_1y_2+x_2y_1\bigr),
\eqno{(\ref{pell.surf.defn}.2)}
$$
which is obtained by identifying  $(x,y)$ 
with  $x+y\sqrt{g} \in k(B)(\sqrt{g})$. 

The {\it inverse} is  $(x,y)\mapsto (x, -y)$; it should be carefully distinguished
from {\it multiplication by $-1$} which is
$(x,y)\mapsto (-x, -y)$. 
\medskip

{\it A 2-valued trivialization \ref{pell.surf.defn}.3.} It is sometimes useful to look at the  2-valued map
$$
T: (x,y,u)\mapsto  \bigl(u, x\pm y\sqrt{g(u)}\bigr)\in B\times \gm, 
$$
where $\gm$ denotes the multiplicative group $\spec k[t,t^{-1}]$.
Its restriction to the $u=b$ fiber is denoted by $T_b$. 
Note that $ x+ y\sqrt{g(u)}$ and $x- y\sqrt{g(u)}$ are inverses.
Thus it makes sense to say that $T_b$ maps a certain point $(x, y,b)$ to a root of unity.
\end{defn}

\begin{defn}[Projective Pell surfaces]\label{proj.pell.surf.defn} 
Continuing with the notation of Definition~\ref{pell.surf.defn}, let  
$S_g\subset \a^2_{xy}\times B$ be an affine  Pell surface and 
$\tilde S_g\subset \p^2_{xyz}\times \bar B$ its closure.
We see in Proposition~\ref{pell.surf.lem.1}  that $\tilde S_g$ is non-normal
along the fiber at infinity  if  $\deg g\geq 2$. We denote its normalization by  $\bar S_g$ and call it the   {\it projective Pell surface}
corresponding to $g$.  The second coordinate projection is
$\pi:\bar S_g\to \bar B$.  

Let $F_{\infty}$ denote the reduced fiber of $\pi$ over $b_{\infty}$ and  
$\bar C_g\subset \bar S_g$  the birational transform of
$\tilde C_g:=(z=0)$.  Thus
$S_g=\bar S_g\setminus\bigl( \bar C+F_{\infty}\bigr)$.

We study the geometry of the pair
$\bigl( \bar S_g, \bar C+F_{\infty}\bigr)$.
These results---especially the computation of the log Kodaira dimension---are not needed for the proofs of the main Theorems, so can be skipped by those who are less interested in the study of open surfaces.
\end{defn}

\begin{defn}[Log Kodaira dimension]\label{log.K.defn}
Let $U$ be a smooth variety and $\bar U\supset U$ a smooth compactification
such that $\bar D:=\bar U\setminus U$ is a divisor with simple normal crossing singularities only.  An {\it $m$-canonical form on $\bar U$ with log poles at infinity} is a section of  $\o_{\bar U}\bigl(mK_{\bar U}+m\bar D\bigr)$ where
$K_{\bar U}$ is the canonical class of $\bar U$. It is easy to see that
the restrictions of  $m$-canonical forms with log poles at infinity to $U$ are independent of the choice of $\bar U$. 

For $m$ sufficiently large and divisible, $P_m(U):=\dim H^0\bigl(\bar U, \o_{\bar U}(mK_{\bar U}+m\bar D)\bigr)$ grows like $(\mbox{constant})m^d$ for some $d\leq \dim U$, called the {\it log Kodaira dimension} of $U$. 
(If $P_m(U)$ is idetically 0, the log Kodaira dimension is declared to be $-1$ by some authors  and $-\infty$ by others.)
We refer to \cite{MR635930} for  basic results on open surfaces  and  to \cite{km-book, kk-singbook} for  a discussions of
 their singularties, including the notion of {\it log canonical pairs.}
\end{defn}

\begin{say}[Special Pell surfaces]\label{special.say}
The Pell equations where $\deg g\leq 1$ are too general and the ones where $g$ is a power of a linear form are too degenerate to be of interest. 
They appear as exceptions to various statements, so we list them here.
\begin{enumerate}
\item $\deg g=0$. Then  $S_g$ is the product of $B$ with a hyperbola. 
\item $\deg g=1$ and $B\cong \a^1$. We denote this by
$S_1:=(x^2-uy^2=1)$. Every other Pell surface is obtained from $S_1$ by pull-back via a
morphism $B\to \a^1$. $S_1$ contains many affine lines, for example
$\bigl(x(t), 1, x(t)^2-1\bigr)$ for any $x(t)\in k[t]$. 
\item $g=ch^2$ is  a constant times a square.  Then 
$x(t)\pm \sqrt{c} h(u(t))y(t)$ are both constants, hence    
so are $x(t)$ and $h(u(t))y(t)$.  Thus the only  (possibly singular) affine lines are
the  obvious ones listed in Paragraph~\ref{aff.lin.aff.var.say}.1.
\item   $g=c(u-\alpha)^d$ for $d\geq 2$.   After a  base field extension we can write these as   $(x^2-u^dy^2=1)$. These  have a $\gm$-action
$(x,y,u)\mapsto  (x, \lambda^{-d} y, \lambda^{2}u)$.
If $d$ is even, this is also a special case of (3). By Example~\ref{all.mult.roots}  the only  (possibly singular) affine lines are
the  obvious ones. 
\end{enumerate}
We see below that the case $\deg g=2$ is also quite special, but
these are very interesting Pell surfaces. 
\begin{enumerate}\setcounter{enumi}{4}
\item  $\deg g=2$  and $B\cong \a^1$.  After a  base field extension we may assume that this is the surface  $S_2=(x^2-(u^2-1)y^2=1)$ that we discussed in Example~\ref{u2-1.exmp}.
\end{enumerate}
\end{say}

\begin{prop} \label{pell.surf.lem.1}  Let $\pi:S_g\to B$ be a Pell surface. 
  The pair $\bigl( \bar S_g, \bar C+F_{\infty}\bigr) $  has the following properties.
\begin{enumerate}
\item $\bar S_g$ is smooth iff $g$ has no multiple roots and $\deg g$ is even.
\item $\bar S_g$ has only $A_m$ singularities.
\item $\bigl( \bar S_g, \bar C+F_{\infty}\bigr) $ is log canonical iff
$g$ has no roots of multiplicity $\geq 3$.
\item $K_{\bar S_g}+\bar  C_g+F_{\infty}\sim  \pi^*\bigl(K_{\bar B}+\tfrac{\deg g+2}{2}[b_{\infty}]\bigr)$.
\item The log Kodaira dimension of $S_g$ is
\begin{enumerate}
\item  $-\infty$  in cases (\ref{special.say}.1--2),
\item $0$  in cases (\ref{special.say}.4--5) and
\item $1$ in all other cases.
\end{enumerate}
\end{enumerate}
\end{prop}

Proof. The fiber of $\pi:\bar S_g\to \bar B$ over $b\in B$ is  smooth
if  $g(b)\neq 0$. Otherwise the fiber  is a pair of lines and, in suitable formal coordinates, a neighborhood of a singular fiber can be written as
$$
\bigl(x^2-z^2-u^my^2=0\bigr)\subset  \p^2_{xyz}\times \hat\a^1_{u},
\eqno{(\ref{pell.surf.lem.1}.6)}
$$
where $\hat\a^1_{u}:=\spec k[[u]]$. 
 If $m\geq 2$ we get an $A_{m-1}$ singularity at $(x=z=u=0)$. 
The curve $\tilde C$ is smooth iff $m=1$ and has an ordinary node iff $m=2$.

Over  $b_{\infty}$ we can rewrite the equation as
$$
\bigl(v^n(x^2-z^2)-y^2=0\bigr)\subset  \p^2_{xyz}\times \hat\a^1_{v}.
$$
This is singular along $v=0$. 
If $n=2r$ is even then  the normalization is given by $y_1=y/v^r$
with equation  
$$
\bigl((x^2-z^2)-y_1^2=0\bigr)\subset  \p^2\times \hat\a^1_{v}.
$$
We have a smooth fiber at $v=0$. If $n=2r+1$ is odd then  the normalization is given by $y_2=y/v^r$
with equation  
$$
\bigl(v(x^2-z^2)-y_2^2=0\bigr)\subset  \p^2\times \hat\a^1_{v}.
$$
We have a double line fiber at $v=0$ and two $A_1$ singular points 
at $v=x\pm z=y_2=0$. These show the claims (1--3). 

The canonical class of $\tilde S_g$ is the restriction of
$\omega_{\p^2}(2)\boxtimes \omega_{\p^1}(n)$. 
At infinity, in the affine chart with equation $v^n(x^2-1)-y^2=0$ a local generator
of the dualizing sheaf is given by $v^{-n}x^{-1}dv\wedge dy$. 
In the $n$ even (resp.\ odd) cases, its pull-back to $\bar S_g$ can be written as
$$
v^{-r}\tfrac{dv\wedge dy_1}{x} \qtq{resp.}
v^{-r}\tfrac{dv\wedge dy_2}{vx} 
$$
Thus $K_{\bar S_g}$ is the pull-back of 
$\omega_{\p^2}(2)\boxtimes \omega_{\p^1}(n-r)$.

 If $n$ is even then  $F_{\infty}\sim \pi^*[b_{\infty}]$ 
and  $n-r=\tfrac{n}{2}$, 
hence
$$
K_{\bar S_g}+\bar  C_g+F_{\infty}\sim  \pi^*\bigl(K_{\bar B}+\tfrac{\deg g+2}{2}[b_{\infty}]\bigr).
$$
 If $n$ is odd then  $2F_{\infty}\sim \pi^*[b_{\infty}]$ 
and  $n-r=\tfrac{n+1}{2}$, 
hence again
$$
2\bigl(K_{\bar S_g}+\bar  C_g+F_{\infty}\bigr)\sim  2\pi^*\bigl(K_{\bar B}+\tfrac{\deg g+2}{2}[b_{\infty}]\bigr).
$$
This proves (4) and also (5) if all roots of $g$ have multiplicity $\leq 2$. 
If $g$ has a root of multiplicity $m_i\geq 3$ then
$\bigl( \bar S_g, \bar C+F_{\infty}\bigr) $ is not log canonical and we need to compute a  log resolution.

Let $c_i$ be the roots of $g$ with multiplicity $m_i\geq 3$ and  
$F_i$ the corresponding fibers. We compute in Claim~\ref{pell.surf.lem.1}.8
that if  $\sigma$ is a section of
$$
\o_{\bar S_g}\bigl(2(K_{\bar S_g}+\bar  C_g+F_{\infty})\bigr)\cong
\pi^*\o_{\bar B}\bigl(2K_{\bar B}+(\deg g+2)[b_{\infty}]\bigr)
$$
that vanishes along $F_i$ with multiplicity  $m_i-2$ then 
 $\sigma$ has only  log poles along $\bar  C_g+F_{\infty}$. 
Since $\sum_i m_i\leq \deg g$, we know that 
 $\sum_i(m_i-2)\leq \deg g-2$, hence such a $\sigma$ exists.
Furthermore, the log Kodaira dimension is $1$ whenever
$\sum_i(m_i-2)< \deg g-2$,
 proving (5) in general. \qed
\medskip

Looking at the last step a little more carefully gives the following more precise version of (\ref{pell.surf.lem.1}.5). 
For $b\in B$ we let $\mult_b(g)$ denote the order of vanishing of $g$ at $p$
and we set $(\mult_b(g)-2)^+=\mult_b(g)-2$ if the latter is positive and $0$ otherwise.
\medskip

{\it Claim \ref{pell.surf.lem.1}.7.}  For $m\geq 1$ 
the $m$-canonical 2-forms on $S_g$ with log poles at infinity are of the form
$$
\begin{array}{l}
\sigma\Bigl(\frac{dy\wedge du}{x}\Bigr)^{\otimes m}\qtq{where}\\
\sigma \in H^0\Bigl(\bar B, \o_{\bar B}\bigl(m\tfrac{\deg g+2}{2}[b_{\infty}]-m\tsum_{p\in B} (\mult_b(g)-2)^+[b]\bigr)\Bigr). 
\end{array}
$$

Proof. We proved that the only non-log-canonical points are the ones in (\ref{pell.surf.lem.1}.6). In affine coordinates we have the pair
$\bigl((x^2-z^2-u^m=0), (z=0)\bigr)$. The  canonical bundle of the surface 
$ (x^2-z^2-u^m=0)$ is generated by  $\tfrac{dx\wedge du}{z} $. Allowing a simple pole along $(z=0)$ and imposing a $c$-fold vanishing along $(u=0)$ gives the form $\sigma_c $ below. We thus need to prove the following.

\medskip

{\it Claim \ref{pell.surf.lem.1}.8.}
Set $U_m:=(x^2-z^2-u^m=0)\subset \a^3$. Then the 2-form
$$
\sigma_c:=u^cz^{-1}\tfrac{dx\wedge du}{z}
$$
has only log poles  iff $c\geq \frac{m-2}{2}$. 
\medskip

Proof. The minimal resolution is covered by charts
$\bigl(x_i=x/u^i, z_i=z/u^i, u\bigr)$ for $i\leq m/2$. The pull-back of
$\sigma_c$ is
$$
u^{c-i}z_i^{-1}\tfrac{dx_i\wedge du}{z_i}.
$$
If $m=2r$ is even then we stop  with $i=r-1$. 
The equation is then $x_{r-1}^2-z_{r-1}^2-u^2=0$  and the origin is a log canonical center of the divisor  $(z_{r-1}=0)$. Thus we need $c\geq r-1$.

If $m=2r+1$ is odd then we stop  with $i=r$. 
The equation is then $x_{r}^2-z_{r}^2-u=0$, hence smooth and, after eliminating $u$, we have the pair
$$
\bigl(\a^2, (z_r=0)+(r-c)(x_{r}^2-z_{r}^2)\bigr).
$$
This is  log canonical  iff  $r-c\leq \frac12$. \qed

\medskip

A surface $T$ is called {\it affine ruled} if there is a dominant morphism
$C\times \a^1\to T$ for some smooth (affine) curve $C$. 
Since the log Kodaira dimension of $\a^1$ is $-\infty$, the
log Kodaira dimension of an affine-ruled surface is also $-\infty$ in characteristic 0, see \cite{MR635930}.
Thus (\ref{pell.surf.lem.1}.5) implies the following. 

\begin{cor} \label{aff.ruled.cor}  Let $k$ be a field of characteristic 0 and 
 $S_g$ a Pell surface over $k$, not isomorphic to one of the special cases (\ref{special.say}.1--2).  Then $S_g$ 
  is not affine-ruled.
In particular, there are only countably many (possibly singular) affine lines on $S_g$. \qed
\end{cor}

\section{Polynomial Pell equations}\label{sec.4}

Much of the theory is already in \cite{abel1826}, a modern treatment with  details and references is in  \cite{scherr} and \cite{zan1, zan2}. 
As in \cite{zan2}, we allow $g(u)$ to have multiple roots, thus the curve $C_g$ defined in 
(\ref{review.ppe}.3) can be singular.

\begin{say}[Review of the theory]\label{review.ppe} 
Let $R$ be an integral domain and 
 $g\in R$. We would like to find  solutions of the (slightly generalized) polynomial Pell equation 
$$
x^2-gy^2=c \qtq{where} x, y\in R,\ c\in R^*.
\eqno{(\ref{review.ppe}.1)}
$$
A lot of the literature on polynomial Pell equations focuses on the case $R=\z[u]$, but here we are interested in the more geometric setting, thus from now on we work over a field $k$ whose characteristic is $\neq 2$
and $R=k[B]$ is the ring of regular functions on a smooth, geometrically irreducible  curve $B$ that has only 1 place at infinity (Definition~\ref{1.pl.at.inf.defn}).  Let $\bar B\supset B$
denote the unique compactification that is smooth at the point at infinity $b_{\infty}$.  
The main example is $B=\a^1$. 

 As we see below, replacing the constant 1 on the right hand side of (\ref{review.ppe}.1) with an arbitrary  $c\in k^*$ is the natural thing to do from the geometric point of view. In the final applications we are mostly interested in algebraically closed fields, and then this does not matter.   
 Note also that if $\phi_1+ \psi_1\sqrt{g} $ is a
 solution of $x^2-gy^2=c$  then
$\frac1{c}\bigl(\phi_1+  \psi_1\sqrt{g} \bigr)^2$ is a
 solution of $x^2-gy^2=1$.

By a solution of (\ref{review.ppe}.1) we mean a pair $\bigl(\phi , \psi \bigr)\in k[B]^2$
for which $\phi ^2-g\psi ^2$ is a nonzero constant. 
Sometimes we call   a function $\phi + \psi\sqrt{g} \in k[B][\sqrt{g}]$  a solution if
$\norm_g (\phi+ \psi\sqrt{g})\in k^*$, where $\norm_g$ denotes the norm of the degree 2 field extension
$k(B)(\sqrt{g})/k(B)$. The advantage of the latter terminology is that all solutions form a multiplicative group. As we see in (\ref{review.ppe}.4), this group is isomorphic to $k^*$ or to $k^*\times \z$.
We aim to describe all solutions up to multiplicative constants. That is, find a generator of  $(k^*\times \z)/k^*\cong \z$.

By looking at the degrees of $\phi ^2$ and of $g\psi ^2$ we see that the only solution is  $\phi=\pm\sqrt{c}, \psi=0$, unless  $\deg g$ is even and  the leading coefficient of $g$ is a square in $k^*$. 
(If $B=\a^1$, we use the usual notion of leading coefficient. Otherwise,
let $v$ be a local parameter at $b_{\infty}$. If $\deg g=n$ then
$v^ng$ is regular and nonzero at  $b_{\infty}$, giving a 
well defined   $(v^ng)(b_{\infty})\in k^*/(k^*)^n$. Thus if $n$ is even then it makes sense to ask whether the leading coefficient of $g$ is a square in $k^*$ or not.)
Note that these hold for $g$ iff they hold after a substitution  $g\bigl(q(t)\bigr)$. In particular, we see the following.
\medskip

{\it Claim \ref{review.ppe}.2.}  If  $\deg g$ is even and    the leading coefficient of $g$ is not a square in $k^*$ then, for every 
nonconstant  $q:B'\to B$, the Pell equation  
$x^2-(g\circ q)y^2=c$ has only obvious solutions (as listed in Paragraph~\ref{aff.lin.aff.var.say}.1). \qed
\medskip

{\it Definition \ref{review.ppe}.3.} 
 Let $C_g$ denote  the  (possibly singular) affine
curve
$C_g:= (v^2=g)\subset \a^1_v\times B$.
Thus  $k[C_g]\cong k[B]+\sqrt{g}k[B]$, as $k[B]$-modules.

 Let  $\bar C_g\supset C_g$ denote its unique projective model  that is smooth at infinity. If $g$ is not a constant times a square then $C_g$ is geometrically irreducible. 
(If $B\cong \a^1$ then $C_g$ is  hyperelliptic.) 

If $\deg g$ is even and the leading coefficient of $g$ is a square in $k$ then  $\bar C_g$  has two $k$-points at infinity; denote these points by $P_1, P_2$. The  following key observation  goes back to 
\cite{abel1826}, but  we
have to pay close attention to the singularities of $C_g$. 
\medskip

{\it Claim \ref{review.ppe}.4.}   
$\bigl(\phi , \psi \bigr)$ is a solution of (\ref{review.ppe}.1) iff  $\phi + \psi\sqrt{g} $ is regular on $C_g$ and  its divisor   is supported on $P_1+P_2$. (See Paragraph~\ref{jac.say} for divisors.)
\medskip

Proof. Note that   $\phi +  \psi\sqrt{g}$ is regular on $C_g$ 
iff $\phi, \psi\in k[B]$. If  $\phi +  \psi\sqrt{g}$ is a solution
then so is
$\phi -  \psi\sqrt{g}$.
Since $\bigl(\phi +  \psi\sqrt{g} \bigr)\bigl(\phi -  \psi \sqrt{g}\bigr)=1$, both factors are units on $C_g$, hence the only possible zeros and poles are at $P_1, P_2$.
Conversely, if  $\phi +  \psi\sqrt{g}$ is regular on $C_g$ and the divisor of $\phi +  \psi \sqrt{g}$ is supported on $P_1+P_2$,
then the same holds for its conjugate. Hence $\phi , \psi $ are both regular functions on $B$  and 
$\phi ^2-g\psi ^2$ is  a regular function  on $B$ without zeros, hence constant, as we noted in Definition~\ref{1.pl.at.inf.defn}. \qed
\medskip

The divisor of $\phi +  \psi \sqrt{g}$ is thus  $m[P_1-P_2]$ for some $m\in \z$. 
This gives an injection
$$
\operatorname{div}\colon  \left\{
\begin{array}{c}
\mbox{solutions of (\ref{review.ppe}.1) up to}\\
\mbox{multiplicative constants}
\end{array}
\right\}\into \z.
\eqno{(\ref{review.ppe}.5)}
$$
(We could also work with $P_2-P_1$,  so the sign involved in
$\operatorname{div}$  is not canonical.)
There is a nontrivial solution iff $[P_1-P_2]\in \jac(\bar C_g)$ is a torsion point.  
(See Paragraph~\ref{jac.say} on Jacobians.)
Its {\it order} is denoted by $\ord(P_1-P_2)$.

A pair  $\bigl(\phi_1, \psi_1\bigr)$ is called a {\it fundamental solution}
iff   $\phi_1,\psi_1$ are regular on $B$ and 
$$
\fdiv\bigl(\phi_1+  \psi_1\sqrt{g}\bigr)=\pm \ord(P_1-P_2)\cdot [P_1-P_2].
\eqno{(\ref{review.ppe}.6)}
$$
(Using the terminology to be  introduced in Definition~\ref{prim.fund.defn}, this holds iff
$\fdiv\bigl(\phi_1+  \psi_1\sqrt{g}\bigr) $ is a  fundamental divisor on $\bar C_g$.)
Up to multiplicative constants,  every other solution is of the form
$$
\pm\phi_n\pm  \psi_n\sqrt{g}= \bigl(\phi_1+  \psi_1\sqrt{g}\bigr)^n
\qtq{for some} n\in \z,
\eqno{(\ref{review.ppe}.7)}
$$
where 
$\bigl(\phi_1+ \sqrt{g} \psi_1\bigr)^{-n}=\bigl(\phi_1- \sqrt{g} \psi_1\bigr)^n$.
Explicitly,
$$
\begin{array}{lcl}
\phi_n &=&\frac{\pm1}{2} \Bigl(\bigl(\phi_1 +  \psi_1\sqrt{g} \bigr)^n+\bigl(\phi_1 -  \psi_1\sqrt{g} \bigr)^n\Bigr),\\[1ex]
\psi_n &=&\frac{\pm1}{2\sqrt{g }} \Bigl(\bigl(\phi_1 +  \psi_1\sqrt{g} \bigr)^n-\bigl(\phi_1 -  \psi_1 \sqrt{g}\bigr)^n\Bigr), \qtq{or}\\[2ex]
\phi_n& = & \pm\sum_{i=0}^{\rdown{n/2}} \binom{n}{2i}\phi_1^{n-2i}\psi_1^{2i}g^i,\\[1ex]
\psi_n& = & \pm\sum_{i=0}^{\rdown{n/2}} \binom{n}{2i+1}\phi_1^{n-2i-1}\psi_1^{2i}g^i.\\
\end{array}
\eqno{(\ref{review.ppe}.8)}
$$
Note that $(\phi_n +  \psi_n\sqrt{g})(\phi_n -  \psi_n\sqrt{g})=1$
implies that  $\phi_n, \psi_n$ and $\phi_n +  \psi_n\sqrt{g}$
have the same order of pole at both $P_1$ and $P_2$.  As noted on \cite{hazama}, this gives the following, 

\medskip

{\it Claim \ref{review.ppe}.9.} Using the above notation, the following hold.   
\begin{enumerate}
\item[(a)] $\deg \phi_n=\deg \psi_n+\tfrac{1}2\deg g=n\cdot \ord(P_1-P_2)$,
\item[(b)]   $\deg \phi_n\geq \tfrac{n}2\deg g$ and
$\deg \psi_n\geq \tfrac{n-1}2\deg g$. 
\item[(c)]  If $(\phi,\psi)$ is a solution and $\deg\psi<\tfrac{1}2\deg g$ then
$(\phi,\psi)$ is a fundamental solution.\qed
\end{enumerate}

\end{say}

We can summarize these considerations as follows.

\begin{cor}  Let $k$ be a field of characteristic $\neq 2$.
Polynomial Pell equations with nontrivial solutions are in  one-to-one correspondence  with pairs $(\bar C\to \bar B, P_1+P_2)$ where
\begin{enumerate}
\item $\bar B$ is a smooth, projective curve over $k$ with a marked point
$b_{\infty}\in \bar B(k)$,
\item $\bar C$ is a reduced, irreducible, projective  curve
equipped with a degree 2 morphism  $\bar C\to \bar B$, 
\item $P_1, P_2\in \bar C(k)$  are the preimages of $b_{\infty}$ and 
\item   $[P_1-P_2]\in \jac(\bar C)$ is a torsion point. 
\end{enumerate}
If these hold then  a rational function  $\Phi$ on $\bar C$ is a
fundamental solution iff    $\Phi$ is regular 
along $\sing \bar C$ and $\fdiv(\Phi)=\pm\ord(P_1-P_2)\cdot [P_1-P_2]$. 
\qed
\end{cor}

  Later we choose Definition~\ref{prim.fund.defn} so that $\Phi$ is a fundamental solution iff its divisor $(\Phi)$ is a fundamental divisor on $\bar C$.

Over a finite field $\f_q$ every point of  $\jac(\bar C)(\f_q)$ is a torsion point, hence, combined with Claim~\ref{review.ppe}.2, we get the following. 

\begin{cor} Let $\f_q$ be  a finite field. The 
Pell equation $x^2-gy^2=1$ has nontrivial solutions
iff  $\deg g$ is even, $g$ is not a constant times a square in $k[B]$ and the leading coefficient of $g$ is a square in $\f_q^*$.\qed
\end{cor}

A geometric constructiom of the correspondence between sections and solutions of the Pell equation is the following.

\begin{prop} \label{pell.=.curves.prop}
Let $\bigl(\bar S, \bar C+F_{\infty}\bigr)$ be a  Pell surface
with identity section $E$ and involution $\tau:\bar C\to \bar C$. Let $\Sigma^+, \Sigma^-$ be   sections given by a nontrivial solution and its  inverse. 
\begin{enumerate}
\item There is a unique rational function $\Phi_S$ with zero along $\Sigma^+$, pole along $\Sigma^-$ and 
value 1 along $E$. 
\item The restriction  $\Phi:=\Phi_S|_{\bar C}$ has zeroes and poles only at
$ \bar C\cap F_{\infty}$.
\item  For $p\in  C$ we have  $\Phi(p)\cdot\Phi\bigl(\tau(p)\bigr)=1$. 
\item The divisor  $(\Phi)$ uniquely determines the pairs $\{\Sigma^+, -\Sigma^+\}$ and $\{\Sigma^-, -\Sigma^-\}$. 
\end{enumerate}
\end{prop}

Proof. 
The formulas for $\Phi_S$ and $\Phi$ are worked out in Paragraph~\ref{conic.compute.say}, they prove (1--3).
The divisor $(\Phi)$ determines $\Phi$ up to a multiplicative constant and
the condition $\Phi(p)\cdot\Phi\bigl(\tau(p)\bigr)=1$ then 
determines $\Phi$ up to sign. It is then again a local computation to show that
$\Phi$ determines $\Phi_S$. The sign ambiguity means that
$\Sigma^+$ and $ -\Sigma^+$ are not distinguished. \qed

\begin{say} \label{conic.compute.say}
Consider the plane conic
$Q:=(x^2-gy^2=z^2)$ and  let $p^+=(a{:}b{:}c)$ and $p^-=(a{:}{-}b{:}c)$ be 
points on it with  $b\neq 0$. Set $e:=(1{:}0{:}1)$. Then 
$$
\Phi_Q:=\frac{(a-c)y-b(x-z)}{(a-c)y+b(x-z)}
$$
is the unique rational function with zero at $p^+$, pole at $p^-$ and 
value 1 at $e$. 

Set $z=0$. Then $x/y=\sqrt{g}$ and the restriction of $\Phi_Q$ to $(z=0)$ becomes
$$
\frac{(a-c)y-bx}{(a-1)y+bx}=\frac{(a-c)-b\sqrt{g}}{(a-c)+b\sqrt{g}}.
$$
Using that $(a-b\sqrt{g})(a+b\sqrt{g})=c^2$, this is further equal to
$$
\frac{a-b\sqrt{g}-c}{c^2(a-b\sqrt{g})^{-1}-c}=-\frac{a}{c}+\frac{b}{c}\sqrt{g}.
$$
In particular,
$$
\Phi_Q(\sqrt{g}{:}1{:}0)\cdot \Phi_Q(-\sqrt{g}{:}1{:}0)=1.
$$
Since a rational function on $Q$ with a single pole is uniquely determined by any 3 of its values, $\Phi_Q(\pm\sqrt{g}{:}1{:}0)$ and $\Phi_Q(1{:}0{:}1)$ also 
determine $\Phi_Q$. 
\end{say}

\begin{say}[Existence of Pell equations with nontrivial solutions]
\label{high.tors.exmp} We give a series of examples where $B\cong \a^1$ and 
$\ord(P_1-P_2)$ is high.

(\ref{high.tors.exmp}.1) If $\bar C$ has genus 1 and $P_1, P_2\in \bar C$ are arbitrary points, the linear system
$|P_1+P_2|$ defines a degree 2 morphism $\bar C\to \p^1$. These correspond to
polynomial Pell equations where $\deg g=4$. We see that the 
order of the torsion can be arbitrary over $\c$.
A complete list of such degree 4 Pell equations over $\q$ is given in
\cite{scherr}; it is quite long and the coefficients of $g$ are complicated.

(\ref{high.tors.exmp}.2) If  $\bar C$ has genus $\geq 2$, then the  hyperelliptic involution $\tau$ is unique. Let 
$\sigma:\bar C\to \jac(\bar C)$ be given by  $P\mapsto  [P-\tau(P)]$. 
Since $\jac(\bar C)$ has countably many torsion points, we
 expect that for very general $\bar C$ the image $\sigma(\bar C)$ does not contain any torsion points. In any case, the image contains   at most finitely many
torsion points 
by \cite{MR717600}.
This was generalized to certain families of curves by \cite{MR3420511}.

Note that  such pairs $(\bar C, P_1+P_2)$ of genus $g$ form a $2g$ dimensional family.
The universal Jacobian over it has dimension $3g$ and the torsion points correspond to a  union of countably many $2g$-dimensional subvarieties. 
Thus,  for every $m$, the family of pairs
$(\bar C, P_1+P_2)$ for which $\ord(P_1-P_2)=m$ 
is either empty or at least  $g$-dimensional.
Next we show that these loci are not empty by constructing higher degree examples where $\ord(P_1-P_2)$  is large.

(\ref{high.tors.exmp}.3) The simplest example is  $x^2-(u^{2m}-1)y^2=1$  with fundamental solution
$u^m-\sqrt{u^{2m}-1}$. Note that
$$
u^m-\sqrt{u^{2m}-1}=u^m\bigl(1-\sqrt{1-u^{-2m}}\bigr)=
u^m\bigl(1-1+\tfrac12 u^{-2m}+\cdots\bigr)=
\tfrac12 u^{-m}+\cdots.
$$
Since $u^{-1}$ is a local parameter at infinity, 
in this exampe $\ord(P_1-P_2)=m$.

(\ref{high.tors.exmp}.4) Assume that  $g_4(u)$ gives $m$ torsion   $(\bar C_4, P_1+P_2)$. 
Corresponding to $g_4(u^n)$  we get  $(\bar C_{4n}, Q_1+Q_2)$,
and the induced map
  $ \bar C_{4n}\to  \bar C_{4}$ is totally ramified at infinity. 
We see in Section~5 that  $\ord(Q_1-Q_2)=nm$, but even the
obvious bound (\ref{3.ques.ques}.1) shows that
$\ord(Q_1-Q_2)$ is a multiple of $m$. 

(\ref{high.tors.exmp}.5) Degree 6 examples might be obtained as follows. Again assume that  $g_4(u)$ gives $m$ torsion.  We may assume that $g_4=ug_3(u)$. 
If  $\phi ^2-ug_3(u)\psi ^2=1$ then 
$x(u^2)^2-g_3(u^2)\bigl(u\psi \bigr)^2=1$. However, the curve $v^2=u^2g_3(u^2)$ is singular, and $v^2=g_3(u^2)$ is smooth. I have not been able to compute how the torsion order changes under normalization of the node.

\end{say}

We can now describe Example~\ref{deg.3.exmp} in terms of the corresponding Jacobians.
 
\begin{say}[Explanation of Example~\ref{deg.3.exmp}]\label{deg.3.exmp.exp}
Start with a Pell surface
$(x^2-g(u)y^2=1)$ where $g(u)$ has odd degree. It has no nontrivial sections by Claim~\ref{review.ppe}.2.

The next simplest thing is to try to find double sections. That is, we look for sections after a degree 2 extension $t=\sqrt{u-c}$. Equivalently,  we use $u=q(t)$ where $q(t)=t^2+c$.  We thus have the diagram
$$
\begin{array}{ccc}
 C_{g\circ q} &\stackrel{\pi_{g\circ q}}{\longrightarrow}& \a^1_t\\
\tau_q\downarrow\hphantom{\tau_q} &&\hphantom{q}\downarrow q\\ 
 C_{g} &\stackrel{\pi_g}{\longrightarrow}& \a^1_u.
\end{array}
$$
Note that $k(C_{g\circ q})=k(u)(\sqrt{u-c}, \sqrt{g})$.
Thus the Galois group of $ k(C_{g\circ q})/k(u)$ is $\z/2\times \z/2$.
Thus there is a 3rd intermediate field
$k(u)(\sqrt{(u-c)g})$ and we have a map
$$
 C_{g\circ q}\to  \bigl(v^2=(u-c)g(u)\bigr)
$$
given by $v\mapsto \sqrt{u-c}\cdot \sqrt{g}$.
Comparing dimensions we see that  $\jac\bigl(v^2=g(t^2+c)\bigr)$   is isogenous to
$$
\jac\bigl(v^2=g(u)\bigr)  \times \jac\bigl(v^2=(u{-}c)g(u)\bigr),
$$
where we use $\jac(\ )$ to denote the Jacobian of the corresponding projective curve that is smooth at infinity.
Thus, although $\jac\bigl(v^2=g(t^2+c)\bigr)$ has dimension $d-1$, 
it is essentially the product of a $\frac{d-1}{2}$-dimensional
Jacobian that is independent of $c$ and of a $\frac{d-1}{2}$-dimensional
Jacobian that varies with $c$. The expectations of Paragraph~\ref{high.tors.exmp} should be 
applied to the family $\bigl(v^2=(u-c)g(u)\bigr)$. 

If $d=3$ then $\bigl(v^2=(u-c)g(u)\bigr)$ is a 1-parameter family of elliptic curves and $P_1-P_2$ is a torsion point for infinitely many values of $c$. 
\end{say}

\begin{say}[Jacobians] \label{jac.say}
We used some facts about Jacobians of singular curves. Many books discuss Jacobians of smooth curves and of nodal curves. However we need to study curves that are geometrically reduced but with worse singularities.
The  exposition  given by Serre \cite{MR0103191} can be easily adapted to our situation. 
The general case is usually treated as a special instance of the theory of Picard varieties outlined in
\cite{FGA}, which seems to be the best reference. See also 
\cite[Chap.9]{blr}.

Let $\bar C$ be a geometrically irreducible and   geometrically reduced curve
over a field $k$  and $C\subset \bar C$ its smooth locus.

By a {\it divisor} on $\bar C$ we mean  a finite linear combination  $D=\sum m_i[c_i]$ where $m_i\in \z$, $c_i\in C$  (not  $\bar C$!).
The degree of a  divisor $D=\sum_i m_i[c_i]$ is 
$$
\deg D:=\tsum_i m_i\deg \bigl(k(c_i)/k\bigr).
$$
Let $f$ be a rational function  on $\bar C$ that is regular and nowhere zero on  $\sing \bar C$. The {\it divisor of $f$} is defined as
$\tsum_{c\in C} v_c(f)[c]$ where $v_c(f)$ is the order of pole
(resp. $-v_c$ is the order of zero)  of $f$ at $c$. 
We do not define the divisor for functions that are either non-regular or vanish at some point of $\sing \bar C$.
The divisor of $f$ is traditionally denoted by $(f)$; we also use $\fdiv(f)$ if confusion is possible.

Two divisors $D_1, D_2$ are  {\it linearly equivalent} if
$D_1-D_2=(f)$ for some  rational function  $f$ on $\bar C$ (that is regular and nowhere zero on  $\sing\bar C$). (Note that while the set of divisors depends only on $C$, the class of rational functions we allow here does depend on the nature of the singularities of $\bar C$.)

The points of $\jac(\bar C)$
are   divisors of  degree 0 on $C$  modulo  linear equivalence. 
The class of a divisor $D$ in the Jacobian is denoted by $[D]$.
We can also think of $\jac(\bar C)$ as parametrizing degree 0 line bundles on $\bar C$.

Later we will need to know that $\jac(\bar C)$ is an algebraic group
of dimension $h^1(\bar C, \o_{\bar C})$. 
If   $\bar C$ is smooth then $ \jac(\bar C)$ is projective, hence   an Abelian variety. Otherwise  $ \jac(\bar C)$ is usually not projective.

We let $\ord(D)$ denote the {\it order} of $D$ as an element of the group
$\jac(\bar C)$.

Note that  our curves $\bar C_g$ are geometrically  irreducible and   geometrically reduced but  singular if $g$ has multiple roots. The arithmetic genus is $\tfrac12 \deg g-1$. 
Thus   $ \jac(\bar C_g)$  has dimension $\tfrac12 \deg g-1$. 

\end{say}

\section{Hazama's treatment of   $x^2-(u^2-1)y^2=1$}\label{sec.5}

First we show in general that finding affine lines on Pell surfaces is
equivalent to describing all sections of some related Pell surfaces.

\begin{say}[Sections and base change] \label{no.new.secs.say}
Let $S_g=(x^2-gy^2=1)$ be a Pell surface over the curve $B$. 
Given any curve $D$, a morphism form $\Phi:D\to S_g$ is given by
a triple  $\Phi=(\phi_x, \phi_y, \phi)$ where
$\phi_x, \phi_y\in k[D]$ and $\phi:D\to B$ is a morphism. 
Alternately, we can view $\Phi$ as a solution of the Pell equation
$x^2-(g\circ \phi) y^2=1$; that is, as  a section of the 
 Pell surface  $S_{g\circ \phi}\to D$. 

Fix now a morphism $\phi:D\to B$
and let   $(x_1(u), y_2(u))$ be a fundamental solution of
$x^2-gy^2=1$. Then   $\bigl(x_1(\phi), y_2(\phi))$ is a  solution of
$x^2-(g\circ \phi)y^2=1$, and so are its powers. If these are all the
solutions of $x^2-(g\circ \phi)y^2=1$ then every lifting
$$
\phi: D\stackrel{\Phi}{\longrightarrow} S_g\to B
\qtq{factors as}
\Phi: D\stackrel{\phi}{\longrightarrow} B\stackrel{(x_n, y_n)}{\longrightarrow} S_g.
$$
Applying this to affine lines shows that the 
 following are equivalent.
\begin{enumerate}
\item On the Pell surface  $S_g=\bigl(x^2-g(u)y^2=1\bigr)$
every (possibly singular) affine line is either vertical or a section.
\item For every nonconstant $q(t)\in k[t]$, 
 $\bigl(\phi_1(q(t)), \psi_1(q(t))\bigr)$ is a 
fundamental solution of  the Pell equation
$x^2-g(q(t))y^2=1$. 
\end{enumerate}
\end{say}

As \cite{hazama} noted, 
this  explains Example~\ref{u2-1.exmp} rather directly. 
His method also proves Theorem~\ref{1.pl.sec.thm}  for this surface.
The only difference is that \cite{hazama} worked  with Jacobians of smooth curves, and these can only handle the cases when  $q(t)^2-1$ has no multiple roots.

\begin{say}[Proof of Example~\ref{u2-1.exmp}] \label{u2-1.exmp.pf}
Let $k$ be a field of characteristic $\neq 2$ and consider the Pell equation
$x^2-(u^2-1)y^2=1$ over $k[u]$.
A  fundamental solution is  $\bigl(\phi , \psi \bigr)=(u, 1)$.
Given any $q(t)\in k[t]$  we get 
the new Pell equation
$$
x^2-\bigl(q(t)^2-1\bigr)y^2=1, 
\eqno{(\ref{u2-1.exmp.pf}.1)}
$$
one of whose solutions is  
$$
\bigl(\phi(q(t)), \psi(q(t))\bigr)=\bigl(q(t), 1\bigr).
\eqno{(\ref{u2-1.exmp.pf}.2)}
$$
 Since the 
fundamental solution is the one with $\deg \psi $ the lowest  (\ref{review.ppe}.9),
it is clear that (\ref{u2-1.exmp.pf}.2) is  a fundamental solution of
(\ref{u2-1.exmp.pf}.1). 

Thus  (\ref{no.new.secs.say}.1--2)   shows that 
every (possibly singular) affine line on $S_2=\bigl(x^2-(u^2-1)y^2=1\bigr)$ is either vertical or a section.
The explicit formula now follows from (\ref{review.ppe}.8). \qed
\end{say}

\section{Torsion order and fundamental index}\label{sec.6}

In this section we study how the 
order and divisibility of divisors changes by pull-back. 
We allow the curves to be singular and the characteristic to be positive.

\begin{defn}\label{prim.fund.defn}
 Let $ \bar C$ be a projective, geometrically reduced and geometrically connected curve over a field $k$. Let $D$ be a Weil divisor supported at smooth points.  $D$ is called {\it primitive} if 
it can not be written as $m'D'$ where $m'>1$ and $D'$ is a Weil divisor.
Thus every Weil divisor $D$ can be uniquely written as
$D=m_1D_1$ where $m_1\geq 1$ and $D_1$ is primitive. We write  $\gcd(D):=m_1$, it is the $\gcd$ of the coefficients of $D$.

We say that $D$ is a {\it principal divisor} if there is  a rational function on $C$ that is regular and  invertible along $\sing  C$, and  such that $(f)=D$.  A principal divisor $D=(f)$---or the function $f$---is called {\it fundamental} if 
it can not be written as $m'D'$ where $m'>1$ and $D'$ is principal.
Thus every principal divisor $D$ can be uniquely written as
$D=m_2D_2$ where  $m_2\geq 1$ and  $D_2$ is fundamental. The value of $m_2$ is called the 
{\it fundamental index} of $D$ or of $f$, and denoted by $\find(D)$ or by $\find(f)$. 
Note that
$$
\gcd(D)=\find(D)\cdot \ord\bigl(D/\gcd(D)\bigr),
\eqno{(\ref{prim.fund.defn}.1)}
$$
where  $\ord(*)$ denotes the  order of $*$ in the group
$\jac(\bar C)$, as in Paragraph~\ref{jac.say}.

{\it Comments.} Neither ``primitive'' nor ``fundamental'' are standard in this context. The notion of primitive coincides with normal usage for vectors in $\z^n$. Fundamental was chosen to coincide with the notion of fundamental solution of Pell's equation, see (\ref{review.ppe}.6).
\end{defn}

In order to prove Theorem~\ref{1.pl.sec.thm}, we need to study the following  question.

\begin{ques}\label{3.ques.ques}
 Let $k$ be a field and  $\pi:\bar C_2\to \bar C_1$  a
flat morphism between  geometrically connected and geometrically reduced curves over $k$. Let  $D_1$  be a divisor on $\bar C_1$ and set
$D_2:=\pi^*D_1$. 
 Assume that  $D_1$ is  fundamental. Is then   $D_2$ also  fundamental?

By (\ref{prim.fund.defn}.1), the answer is related to the change of the order under pull-back. 
Pulling pack  of a rational function on $\bar C_1$ and
taking the norm of a rational function on $\bar C_2$ shows that
$$
\ord(D_2)\mid  \ord(D_1)\qtq{and}\ord(D_1)\mid \deg \pi \cdot\ord(D_2).
\eqno{(\ref{3.ques.ques}.1)}
$$
 However, we need more precise information. 
\end{ques}

The following lemma connects the fundamental index to \'etale covers.

\begin{lem}\label{find.equiv.covers}
 Let $\bar C$ be a  projective, geometrically reduced and geometrically connected curve over a field $k$. Let $f$ be a rational function that is
regular and  invertible along $\sing \bar C$. Set $D:=(f)$ and $C:=\bar C\setminus D$. 
Fix $m\in \n$ not divisible by $\chr k$ and set 
$$
 C[\sqrt[m]{f}]
:= \bigl( u^m=f\bigr)\subset  C\times \a^1_{u}.
\eqno{(\ref{find.equiv.covers}.1)}
$$
The following are equivalent.
\begin{enumerate}\setcounter{enumi}{1}
\item There is a   regular  function $g$ on $C$ such that $g^m=cf$ for some $c\in k^*$.
\item $m\mid \find(f)$.
\item The projection $ C[\sqrt[m]{cf}]\to C$ has a section for some $c\in k^*$.
\item The projection $ C[\sqrt[m]{cf}]\to C$ has a section for some $c\in k^*$
that is also a connected component of $ C[\sqrt[m]{cf}]$. 
\end{enumerate}
\end{lem}

Proof. If $g^m=cf$ then $(f)=m(g)$, hence (2) $\Rightarrow$ (3).
Conversely, if $(g)=\frac1{m}(f)$ then  $g^m=cf$ for some $c\in k^*$.

If $g^m=cf$ then $(u=g)$ defines a section and if 
 $\sigma: C\to  C[\sqrt[m]{cf}]$ is a section then $g:=u\circ \sigma$ satisfies
 $g^m=cf$ for some $c\in k^*$.
Since $m$ is not divisible by $\chr k$, the projection $ C[\sqrt[m]{cf}]\to C$ is \'etale, so any section is also a connected component. \qed
\medskip

We can now prove the following criterion for the preservation of the
fundamental index by pull-backs in characteristic 0.

\begin{thm} \label{H1.surj.fund.sols.thm}
Let $ \bar g:\bar C_2\to \bar C_1$ be a finite morphism of 
projective,  reduced and connected curves over  $\c$. Let $\phi_1$ be a rational function on $\bar C_1$ and
$\phi_2:=\phi_1\circ \bar g$. Assume that $\phi_i$ is regular and invertible along $\sing \bar C_i$ for $i=1,2$. Set $C_i:=\bar C_i\setminus \supp (\phi_i)$; by restricion we get a finite morphism $g:C_2\to  C_1$. 

 If $g_*:H_1(C_2,\z)\to  H_1(C_1,\z)$ is surjective then
$\find(\phi_1)=\find(\phi_2)$.
\end{thm}

Proof. We may as well assume that $\find(\phi_1)=1$.
Pick a prime $\ell$  and,
as in Lemma~\ref{find.equiv.covers}, consider the cover
$ C_1[\sqrt[\ell]{\phi_1}]\to C_1$. It is
\begin{enumerate}
\item connected by 
(\ref{find.equiv.covers}.5) $\Rightarrow$ (\ref{find.equiv.covers}.3),
\item  Galois cover with Galois group
$\mu_\ell$, the group of $\ell$th roots of unity acting by multiplication and
\item unramified
since $\phi_1$ has neither zeros nor poles on $C_1$.
\end{enumerate}
These imply that  $ C_1[\sqrt[\ell]{\phi_1}]\to C_1$  corresponds to a surjective homomorphism
$\sigma_1: H_1(C_1,\z)\to \mu_\ell$.
The map of the fiber product
$$
C_2[\sqrt[\ell]{\phi_2}]\cong C_2\times_{C_1}C_1[\sqrt[\ell]{\phi_1}]\to C_2
$$
then corresponds to the composite
$$
\sigma_2: H_1(C_2,\z)\stackrel{g_*}{\to} H_1(C_1,\z)\stackrel{\sigma_1}{\to} \mu_\ell,
$$
which is also surjective if $g_*$ is surjective. 
Therefore $C_2[\sqrt[\ell]{\phi_2}]$
is also connected.   
Thus $\find(\phi_2)$ is not divisible by $\ell$ by
(\ref{find.equiv.covers}.3) $\Rightarrow$ (\ref{find.equiv.covers}.5). 
We conclude by using this for every $\ell$. \qed

\medskip

It is not hard to prove a version of this in positive characteristic, at least for primes other than the characteristic. However, first we concentrate on characteristic 0 and, in the next section, we aim 
 to understand when maps between algebraic varieties induce a  surjection on the first homology groups.

\section{$H_1$-surjective maps}\label{sec.7}

We start with the simplest statement that is needed for 
the proof of Theorem~\ref{1.pl.sec.thm} over $\c$.

\begin{prop}\label{pi1.onto.curve.prop} Consider the fiber product diagram
$$
\begin{array}{ccc}
 D\times_{\c}B & \stackrel{g_B}{\longrightarrow} & B\\
h_D\downarrow\hphantom{h_D}  &&  \hphantom{h}\downarrow h\\
 D  & \stackrel{g}{\longrightarrow} &  \c,
\end{array}
$$
where $B$ is a smooth, connected curve, $D$ is a connected,
possibly singular curve and 
  $g:D\to \c$ and $h:B\to \c$ are proper morphisms. Then  $D\times_{\c}B$ 
is connected and the induced map
$\pi_1\bigl(D\times_{\c}B\bigr)\to \pi_1(D)$  is surjective. Thus
  $$
 H_1\bigl(D\times_{\c}B, \z\bigr)\to H_1(D,\z)
  \qtq{is also surjective.}
  $$
\end{prop}

{\it Note.} We need surjectivity for $H_1$, so proving
surjectivity for $\pi_1$ seems overkill. However, we see in Example~\ref{com.H1.pi1.3} that
  one has to focus on  $\pi_1$ and switch to $H_1$ only at the very end of the proof. 
\medskip

We start by establishing some topological properties of the maps
 $g:D\to \c$ and $h:B\to \c$, and then prove
Proposition~\ref{pi1.onto.curve.prop} using only these.
We start with $g:D\to \c$.

\begin{say}\label{weak.l.p.say}
  We say that a continuous map of topological spaces $g:M\to N$ has the {\it 
path lifting property} if the following holds.
\begin{enumerate}
\item  Given any continuous map  $\gamma:[0,1]\to N$ and $m\in M$ such that 
$g(m)=\gamma(0)$,
  there is a continuous map  $\gamma':[0,1]\to M$ such that $\gamma'(0)=m$ and 
$\gamma=g\circ \gamma'$. We do not require $\gamma'$ to be unique.
\end{enumerate}
  Every proper, flat surjection of (possibly singular) curves over $\c$ 
has the  path lifting property. Thus our map $g:D\to \c$ in Proposition~\ref{pi1.onto.curve.prop} has the  path lifting property.
More generally, every proper, universally
open, pure relative dimensional, surjective  morphism of
  $\c$-schemes of finite type has the  path lifting property, see  \cite[Sec.3]{k-pi1}.

  Assume that $g:M\to N$ has the  path lifting property and, in addition, $g^{-1}(n)$ is finite for some $n\in N$.
  Pick $m\in g^{-1}(n)$.
  Then every loop $\gamma$ starting and ending in $n$ lifts to a path that 
starts at $m$ and ends in $g^{-1}(n)$. If 2 loops  $\gamma_1, \gamma_2$ end
at the same point then $\gamma_1\gamma_2^{-1}$ lifts to a loop on $M$. This
shows that
  the image of $\pi_1(M,m)\to \pi_1(N,n) $ has finite index in $\pi_1(N,n) 
$.
  \end{say}

\begin{lem} \label{triv.paths.lem}
Let  $B$ be a smooth, connected curve and 
  $h:B\to \c$ a proper morphism.   Pick general $c\in \c$ with preimages $m_1,\dots, m_d$. Then, for every $1\leq i,j\leq d$ there are
paths  $\phi_{ij}, \psi_{ij}$ such that
$$
\phi_{ij}(0)=m_i, \phi_{ij}(1)=m_j,\ \psi_{ij}(0)=m_j, \psi_{ij}(1)=m_j
\qtq{and} h\circ \phi_{ij}=h\circ \psi_{ij}.
$$
\end{lem}

Proof. 
Let $p_1,\dots, p_r\in \c$ be the branch points of $h$. We may as well assume that the line segments  $[p_k, c)$ are disjoint.
Then $h^{-1}[p_k, c]$ is a union of paths  $\gamma_{k\ell}:[0,1]\to B$
such that  $h\circ \gamma_{k\ell}(0)=p_k$, $g\circ \gamma_{k\ell}(1)=c$ and 
$h\circ \gamma_{k\ell}$ is independent of $\ell$. 
Construct a   graph with vertices $m_1,\dots, m_d$ where 2 vertices $m_i, m_j$  are connected by an edge if there is a $k$ and $\ell_i, \ell_j$ such that
$$
  \gamma_{k\ell_i}(0)=\gamma_{k\ell_j}(0),\
\gamma_{k\ell_i}(1)=m_i\qtq{and} \gamma_{k\ell_j}(1)=m_j.
\eqno{(\ref{triv.paths.lem}.1)}
$$
Thus  $\phi_{ij}=\gamma_{k\ell_i}\ast \gamma_{k\ell_j}^{-1}$ and
$\psi_{ij}=\gamma_{k\ell_i}\ast \gamma_{k\ell_j}^{-1}$ work for the pair $m_i, m_j$,
where $\ast$ denotes the concatenation of arcs.
The graph is connected
since $B$ is connected. Thus suitable concatenations of  the above pairs give a solution for every
$m_i, m_j$. \qed

\begin{say}[Proof of Proposition~\ref{pi1.onto.curve.prop}]
\label{pi1.onto.curve.prop.pf}
    Pick general $c\in \c$ with preimages $m_1,\dots, m_d\in B$. Choose
$m_i, m_j$ and let   $\phi_{ij}, \psi_{ij}$ be as in
(\ref{triv.paths.lem}).

Pick any  $n\in g^{-1}(c)$ and
let $\gamma'$ be a lifting of $h\circ \phi_{ij}=h\circ \psi_{ij}$ to $D$, going
from  $n$ to another point $n'$.
Set $\phi'_{ij}:=(\gamma',\phi_{ij})$ and $\psi'_{ij}:=(\gamma', \psi_{ij})$.
Note that
$$
\phi'_{ij}(0)=(n,m_i),\ \phi'_{ij}(1)=(n',m_j),\ \psi'_{ij}(0)=(n,m_j),\
\psi_{ij}(1)=(n',m_j).
$$
The concatenation of $\phi'_{ij}$ with the inverse of $\psi'_{ij}$
is a path in $D\times_{\c}B$ that starts at $(n,m_i)$ and ends at $(n,m_j)$.
Thus $D\times_{\c}B$ is connected.

Thus, as we noted in Paragraph~\ref{weak.l.p.say},
the image of $\pi_1\bigl(D\times_{\c}B\bigr)\to \pi_1(D)$ has finite index
in
$\pi_1(D)$. Let $D'\to D$ be the corresponding covering space.

We can apply the above argument to $h:B\to \c$ and $g':D'\to \c$
to conclude that $D'\times_{\c}B$ is connected.  On the other hand,
the number of its connected components is the index of
 $\im\bigl[\pi_1\bigl(D\times_{\c}B\bigr)\to \pi_1(D)\bigr]$  in
$\pi_1(D)$. Thus $\pi_1\bigl(D\times_{\c}B\bigr)\to \pi_1(D)$ is surjective.
\qed
\end{say}

\section{Proof of Theorem~\ref{1.pl.sec.thm} over $\c$}
\label{sec.8}

We start over any field and then we point out where the
characteristic 0 assumption is used.

\begin{say}[Proof of Theorem~\ref{1.pl.sec.thm}]\label{basic.setup.say}
Start with a  Pell equation
$x^2-g(u)y^2=1$. We may as well assume that $k$ is  algebraically closed.
As we noted in Paragraph~\ref{no.new.secs.say},
it is enough to show that  for every proper morphism $q:B\to \a^1$,  all
solutions of the new Pell equation
$x^2-(g\circ q)y^2=1$  come from a solution of $x^2-g(u)y^2=1$.

Now we switch to the geometric side described in  Proposition~\ref{pell.=.curves.prop}.

As in  (\ref{review.ppe}.3)
let $\bar C_g$  and $\bar C_{g\circ q}$ be the corresponding curves
and $P_1, P_2\in \bar C_g$ and  $Q^{(q)}_1, Q^{(q)}_2\in \bar C_{g\circ q}$ the points at infinity. We get a commutative diagram
$$
\begin{array}{ccccc}
Q^{(q)}_1+ Q^{(q)}_2 & \subset & \bar C_{g\circ q} &\stackrel{\pi_{g\circ q}}{\longrightarrow}& \bar B\\
\downarrow &&\hphantom{\tau_q}\downarrow\tau_q &&\hphantom{q}\downarrow q\\ 
P_1+P_2 & \subset & \bar C_{g} &\stackrel{\pi_g}{\longrightarrow}& \p^1_u
\end{array}
\eqno{(\ref{basic.setup.say}.1)}
$$
where the right hand side is a fiber product square.
Since $\deg g$ is even,   $\pi_g$ and  $\pi_{g\circ q}$ are \'etale over the points at infinity. Since $q$ has  ramification index $=\deg q$ at infinity, 
 $\tau_q$ also has  ramification index $=\deg q$ at $Q^{(q)}_1, Q^{(q)}_2$.

We distinguish 2 cases,  depending on $\ord(P_1- P_2)$. 

\medskip

{\it Non-torsion case \ref{basic.setup.say}.2.} If 
$\ord(P_1- P_2)=\infty$ then also  $\ord(Q^{(q)}_1- Q^{(q)}_2)=\infty$ by (\ref{3.ques.ques}.1), hence
$x^2-(g\circ q)y^2=1$  has only trivial solutions by
  Proposition~\ref{pell.=.curves.prop}.  Thus Theorem~\ref{1.pl.sec.thm}
holds in this case.

\medskip

{\it Torsion case \ref{basic.setup.say}.3.}
 If  $\ord(P_1- P_2)=n$ is finite  then $x^2-g(u)y^2=1$ has 
nontrivial solutions by Proposition~\ref{pell.=.curves.prop}.
Thus it has 
a fundamental solution $$\Phi:=x_1(u)+ y_1(u)\sqrt{g(u)},
$$
whose divisor is  $D_g:=(\Phi)=n(P_1- P_2)$.  
As we noted in Paragraph~\ref{no.new.secs.say},
it is sufficient to show that 
$$
\Phi_q:=x_1(q(t))+ y_1(q(t))\sqrt{(g\circ q)(t)}
$$ 
is  a  fundamental solution of $x^2-(g\circ q)y^2=1$ for every $q$.

{\it Step \ref{basic.setup.say}.4.} If we are over $\c$, then, 
by Theorem~\ref{H1.surj.fund.sols.thm}, $\Phi_q$  is  a  fundamental solution if
$\tau_q: C_{g\circ q}\to C_g$ induces a surjection on the first (topological) homology groups.  We have a fiber product diagram
$$
\begin{array}{ccc}
 C_{g\circ q} &\stackrel{\pi_{g\circ q}}{\longrightarrow}& B\\
\tau_q\downarrow\hphantom{\tau_q} &&\hphantom{q}\downarrow q\\ 
 C_{g} &\stackrel{\pi_g}{\longrightarrow}& \a^1_u.
\end{array}
$$
By Proposition~\ref{pi1.onto.curve.prop}, $\tau_q: C_{g\circ q}\to C_g$ induces a surjection on the
fundamental groups. Since the first homology group is the abelianization of the fundamental group, we see that
$$
(\tau_q)_*: H_1\bigl(C_{g\circ q},\z\bigr)\to H_1\bigl(C_g,\z\bigr)
\qtq{is also surjective.}
$$ 
This completes the proof of Theorem~\ref{1.pl.sec.thm} in characteristic 0. \qed
\end{say}

\section{Theorem~\ref{1.pl.sec.thm} in positive characteristic}
\label{sec.9}

\begin{say}[Proof of Theorem~\ref{1.pl.sec.thm}]\label{basic.setup.say.p}
In positive characteristic, we start the proof exactly as in
Paragraph~\ref{basic.setup.say}. Everything works as before until we reach Step~\ref{basic.setup.say}.4.
At this point we have the fiber product diagram
$$
\begin{array}{ccc}
 C_{g\circ q} &\stackrel{\pi_{g\circ q}}{\longrightarrow}& B\\
\tau_q\downarrow\hphantom{\tau_q} &&\hphantom{q}\downarrow q\\ 
 C_{g} &\stackrel{\pi_g}{\longrightarrow}& \a^1_u,
\end{array}
$$
and we would like to prove that 
$$
(\pi_{g\circ q})_*: H_1(C_{g\circ q})\to H_1(C_{g})\qtq{is surjective,}
$$
where we define the {\it algebraic first homology group} $H_1(X)$ of a scheme $X$ as the abelianization of the algebraic fundamental group  $\pi_1(X)$.

A new problem we face is that the while the projective line $\p^1$ is simply connected, the affine line $\a^1$  is not simply connected in positive characteristic.  In fact, $\pi_1(\a^1_k)$ is a very large group which depends on $k$ and it has not been fully determined. 

We have  a rather complicated way of getting around this issue.

\medskip
{\it Step \ref{basic.setup.say.p}.1.} As a direct analog of
Proposition~\ref{pi1.onto.curve.prop} we show that
if   $q_*:\pi_1(B)\to \pi_1(\a^1)$  is surjective 
then so is 
$$
(\pi_{g\circ q})_*: \pi_1(C_{g\circ q})\to \pi_1(C_{g}).
$$
This turns out to be a rather general property of certain fiber product digrams;
see Section~\ref{pi1.surj.sec} for details. 

\medskip
{\it Step \ref{basic.setup.say.p}.2.}
While not every map $q:B\to \a^1$ is $\pi_1$-surjective,
we show in Section~\ref{pi1.surj.crit.sect}
that all maps with sufficiently `mild' ramification at infinity are $\pi_1$-surjective.  This proves  Theorem~\ref{1.pl.sec.thm} whenever $q:B\to \a^1$
has `mild' ramification at infinity. 

\medskip
{\it Step \ref{basic.setup.say.p}.3.}
We show that a general deformation of any $(q:B\to \a^1)$ 
has `mild' ramification at infinity, hence it is 
is $\pi_1$-surjective. This is rather basic deformation theory;
see Section~\ref{def.morph.sec}.

\medskip
{\it Step \ref{basic.setup.say.p}.4.}
We prove in  Section~\ref{fund.ind.def.sec} that the prime-to-$p$ part of
$\find \bigl(\pi^*_{g\circ q}(P_1-P_2)\bigr)$ is unchanged by deformations. 
Combining this with Steps~\ref{basic.setup.say.p}.2--3  we obtain that 
the pull-back of a fundamental solution never becomes an $m$th power 
for $p\nmid m$. 

\medskip
{\it Step \ref{basic.setup.say.p}.5.} It remains to show that 
the pull-back of a fundamental solution never becomes a $p$th power.  
This follows from Proposition~\ref{p.power.no.prop}. 
Note that inseparable multisections have been especially troublesome for elliptic K3 surfaces; see \cite{2019arXiv190404803B} for a discussion.
\qed
\end{say}

\begin{prop}\label{p.power.no.prop}
 Let $k$ be a perfect field of odd characteristc $p$, $B$ a smooth curve with 1 place at infinity and   $q:B\to \a^1$  a  finite morphism.  
Let  $x_1+y_1\sqrt{g}$  be a  solution of a Pell equation
$x^2-g(u)y^2=1$ in $k[u]$. If  $x_1+y_1\sqrt{g}$ is  a $p$th power in 
 $k[B][\sqrt{g}]$ then it is also a $p$th power in 
 $k[u, \sqrt{g}]$.
\end{prop}

{\it Warning.} If $q$ is purely inseparable then  $x_1+y_1\sqrt{g}$ is  always a $p$th power in the function field
 $k(B)[\sqrt{g}]$.  However the ring $k[B][\sqrt{g}]$ is not normal, so this does not contradict our claim.
\medskip

Proof.  We use induction on the degree of inseparability of $q$.

If $q$ is separable, then so is  $k(B, \sqrt{g})/k(u, \sqrt{g})$,
so if an $h\in k(u, \sqrt{g})$ is a $p$th power in   $k(B, \sqrt{g})$
then it is already a $p$th power in   $k(u, \sqrt{g})$.
That is, $x_1+y_1\sqrt{g}=\bigl(x_0+y_0\sqrt{g}\bigr)^p$ where
$x_0,y_0\in k(u)$. On the other hand, $ x_0+y_0\sqrt{g}$ is also the only possible $p$th root in $k(B, \sqrt{g})$, so
$x_0,y_0\in k[B]$.  
Thus $x_0,y_0$ are integral over $k[u]$.
Since $k[u]$ is integrally closed in $k(u)$, we get that   $x_0,y_0\in k[u]$.

If $q$ is not separable, we can factor it as
$$
q: B\stackrel{q'}{\longrightarrow} \a^1\stackrel{F}{\longrightarrow}\a^1,
$$
where $F$ is  the Frobenius.
Since the degree of inseparability of $q'$ is less than the degree of inseparability of $q$, by induction   $x_1+y_1\sqrt{g}=\bigl(x_0+y_0\sqrt{g}\bigr)^p$ where
$x_0,y_0\in k[u^{1/p}]$.
 Taking $p$th powers, we get that
$$
(x_0^p)^2-g(u)^p(y_0^p)^2=1\qtq{and} x_0^p, y_0^p\in k[u].
$$
Lemma~\ref{all.mult.roots.cor} now gives 
that $x_0,y_0\in k[u]$. \qed

\begin{lem}\label{all.mult.roots.cor}
Let $k$ be a perfect field of characteristic $p> 2$. 
Then every solution of 
 $x^2-g^p(u)y^2=1$  in $k[u]$ is of the form
$(x_2^p, y_2^p)$ where  $x_2^2-g(u)y_2^2=1$. 
\end{lem}

Proof. Every root of $g^p$ is multiple, so $x(u)$ is a $p$th power by
Example~\ref{all.mult.roots} and then so is $y(u)$. \qed

\section{$\pi_1$-surjective maps}\label{pi1.surj.sec}

In Step~\ref{basic.setup.say}.4 of the proof of Theorem~\ref{1.pl.sec.thm}
 it would be useful to know that
being surjective on the first  homology group is preserved by base change.
This is, however, not true, see Example~\ref{com.H1.pi1.3}.
By contrast, we get much better behaviour for the fundamental group,
as shown by the next result of \cite[Sec.1]{k-pi1}.

\begin{thm} \label{com.H1.pi1.1} Let $k$  be a field and consider a fiber product diagram
$$
\begin{array}{ccc}
 X\times_SY & \stackrel{g_Y}{\longrightarrow} & Y\\
\downarrow  &&  \hphantom{h}\downarrow h\\
 X  & \stackrel{g}{\longrightarrow} &  S,
\end{array}
\eqno{(\ref{com.H1.pi1.1}.1)}
$$
where  $X, Y, S$   are  geometrically connected $k$-schemes
and $g, h$ are finite,  universally open morphisms. 
Assume that  $g$  induces
a surjection on the fundamental groups.

Then 
$X\times_SY$ is  geometrically connected and 
 $g_Y$ also induces
a surjection on the fundamental groups. \qed
\end{thm}

\begin{exmp}\label{com.H1.pi1.3}  
Let $X$ be a simply connected manifold (or variety over $\c$)  on which the alternating group $A_n$ acts freely. Assume that $n\geq 6$ is odd. 
Let $A_{n-1}\subset A_n$ be a point stabilizer and  $C_n\subset A_n$  a subgroup generated by an $n$-cycle. We get a commutative diagram
$$
\begin{array}{ccc}
X & \stackrel{g'}{\longrightarrow} & X/C_n\\
\downarrow && \downarrow \\
X/A_{n-1} & \stackrel{g}{\longrightarrow} & X/A_n,
\end{array}
\eqno{(\ref{com.H1.pi1.3}.1)}
$$
which is a fiber product square.  Since $n\geq 6$, 
$A_{n-1}$ and $A_n$ are simple, so 
$H_1\bigl(X/A_{n-1},\z\bigr)$ and $H_1\bigl(X/A_{n},\z\bigr)$ are both trivial.
Thus $g$ is  $H_1$-surjective. 
However
$$
g'_*:  H_1\bigl(X,\z\bigr)\to H_1\bigl(X/C_n,\z\bigr)\cong C_n
$$
is not  surjective since $ H_1(X,\z)$ is trivial.
\end{exmp}
 \medskip

\begin{exmp}\label{com.H1.pi1.4}   If $k$ is algebraically closed and $0<\chr k<n$ then
then $\a^1_k$ has \'etale, Galois covers with Galois group $A_n$ by
\cite{Raynaud1994}. 
We can thus obtain a base change diagram as (\ref{com.H1.pi1.3}.1)
(though not with $X$ simply connected).  If $p$ does not divide $n$
then we get a diagram 
$$
\begin{array}{ccc}
A\times_{\a^1}B & \stackrel{g'}{\longrightarrow} & B\\
\downarrow && \downarrow \\
A & \stackrel{g}{\longrightarrow} & \hphantom{x}\a^1,
\end{array}
\eqno{(\ref{com.H1.pi1.4}.1)}
$$
where $g$ is surjective on the algebraic $H_1$  (up to $p^{\infty}$-torsion)
but $g'$ is not. 
\end{exmp}

\section{Criterion for $\pi_1$-surjectivty}\label{pi1.surj.crit.sect}

We prove a  condition of 
$\pi_1$-surjectivity, in terms of the discriminant at infinity.

\begin{say}[Discriminant]\label{discr.say}
 Let $g:C\to B$ be a separable morphism between smooth,
 projective curves over a field $k$ of characteristic $p\geq 0$. The sheaf $\omega_C/g^*\omega_B$ or---more frequently---its associated divisor
$$
\mfrd(g):=\tsum_c \mfrd_c(g)[c]:=\tsum_c \dim_{k(c)}(\omega_C/g^*\omega_B)[c]
\eqno{(\ref{discr.say}.1)}
$$
is called the {\it discriminant} of $g$.
Thus
$$
\begin{array}{rcl}
\deg \mfrd(g)&=&\deg \omega_C-\deg g^* \omega_B\\
&=&\deg \omega_C-\deg g\cdot\deg \omega_B.
\end{array}
\eqno{(\ref{discr.say}.2)}
$$
Pick points $c\in C$, $b=g(c)$ and local coordinates  $s$ at $c$ and $t$ at $b$. We can then write
$g^*t=\phi(s)$ for some  function $\phi$ that is regular and vanishes at $c$. 
The {\it ramification index} of $g$ at  $c$ is  $e_c(g):=\mult_c \phi(s)$.

Since $g^*dt=d\bigl(\phi(s)\bigr)=\phi'(s)ds$, we see that
$$
\mfrd_c(g)=\mult_c \phi'(s).
\eqno{(\ref{discr.say}.3)}
$$
This shows that
$$
\mfrd_c(g)\geq e_c(g)-1
\qtq{and equality holds iff} p\nmid e_c(g).
\eqno{(\ref{discr.say}.4)}
$$
We say that $g$ is {\it tamely ramified} at $c$ if $\mfrd_c(g)=e_c(g)-1$
and {\it wildly ramified} at $c$ if $\mfrd_c(g)>e_c(g)-1$.
Note that  $g$ is wildly ramified at $c$ iff  $p\mid e_c(g)$.

Let $g_i:(C_i, c_i)\to (C_{i+1}, c_{i+1})$ be morphisms of smooth, pointed curves.  Choose local coordinates $s_i$ at $c_i$. Then $g_i$ can be given
as  $g_i^*s_{i+1}=\phi_i(s_i)$. Thus
$(g_2\circ g_1)^*s_3=\phi_2\bigl(\phi_1(s_1)\bigr)$ and so
$$
(g_2\circ g_1)^*ds_3=\phi'_2\bigl(\phi_1(s_1)\bigr)\cdot \phi'_1(s_1) \cdot ds_1.
\eqno{(\ref{discr.say}.5)}
$$
Taking the multiplicity at $c_1$ gives the formula
$$
\mfrd_{c_1}(g_2\circ g_1)=\mfrd_{c_2}(g_2)e_{c_1}(g_1)+\mfrd_{c_1}(g_1).
\eqno{(\ref{discr.say}.6)}
$$
\end{say}

\begin{lem}\label{pi1.onto.crit.lem}
 Let $C$ be a smooth projective curve over a  field of characteristic
 $p>0$ and 
$g:C\to \p^1$   a separable morphism
such that $g^{-1}(\infty)=\{c\}$ is a single point. Assume that
$$
\mfrd_c(g)<2\bigl(1-\tfrac1{p}\bigr) \deg g.
\eqno{(\ref{pi1.onto.crit.lem}.1)}
$$
Then $g_*:\pi_1(C\setminus\{c\})\to \pi_1(\a^1)$ is surjective.
\end{lem}

Proof. If $g_*:\pi_1(C\setminus\{c\})\to \pi_1(\a^1)$ is not surjective
then $g$ factors as
$$
g:  (C, c)\stackrel{r_1}{\to} (B, b)\stackrel{r_2}{\to}  (\p^1, \infty)
$$
where $r_2: B\setminus\{b\}\to \a^1$ is \'etale and $\deg r_2\geq 2$. By the Hurwitz formula
$$
\mfrd_b(r_2)=2\deg r_2 +2g(B)-2\geq 2\deg r_2 -2=
2\bigl(1-\tfrac1{\deg r_2}\bigr) \deg r_2.
$$
In paticular, $ \mfrd_b(r_2) \geq \deg r_2$ and so $r_2$ has wild ramification 
at $b$. Therefore  $\deg r_2$  is divisible by $p$. Thus we obtain that
$$
\mfrd_b(r_2)\geq 2\bigl(1-\tfrac1{p}\bigr) \deg r_2.
$$
Combining this with (\ref{discr.say}.6) we get that
$$
\begin{array}{rcl}
\mfrd_c(g)&\geq& 2\bigl(1-\tfrac1{p}\bigr) \deg r_2\deg r_1+\mfrd_c(r_1)\\[1ex]
&\geq& 2\bigl(1-\tfrac1{p}\bigr) \deg g. \qed
\end{array}
$$

\section{Deformation of morphisms to $\a^1$}\label{def.morph.sec}

In positive characteristic we still need to deal with morphisms
$B\to \a^1$ that are not $\pi_1$-surjective. The next result says that
a suitable small deformation of any $B\to \a^1$ is $\pi_1$-surjective.

\begin{defn} The {\it Hurwitz scheme} ${\mathcal H}_{d,g}$ paramerizes 
degree $d$ morphisms  $C\to \p^1$ from a smooth, projective curve of genus $g$ to $\p^1$; see \cite{MR0260752, MR1837111}.

Let ${\mathcal H}_{d,g;d}\subset {\mathcal H}_{d,g}$ denote the closed subset
parametrizing those maps $\bar\pi:\bar B\to \p^1$ for which
$\bar\pi^{-1}(\infty)$ consists of a unique point, denoted by $b_{\infty}$.
The  ramification index of $\bar\pi$ equals $d$ at $b_{\infty}$.

Thus $B:=\bar B\setminus\{b_{\infty}\}$ is a smooth curve of genus $g$
with 1 place at infinity and $\pi:=\bar \pi|_B:B\to \a^1$
is a finite morphisms of degree $d$. Thus 
${\mathcal H}_{d,g;d}$ is also the moduli space of 
 genus $g$ curves
with 1 place at infinity, equipped with a finite morphism
$\pi:B\to \a^1$ of degree $d$.

Over $\c$ the
Hurwitz scheme ${\mathcal H}_{d,g}$ is irreducible; historically this gave the first proof that the moduli space of genus $g$ curves is irreducible \cite{MR1510692}.
See  \cite{MR828827} for a purely topological approach.

It is natural to hope that ${\mathcal H}_{d,g;d}$ is irreducible over any field.  A positive answer would give a shorter  proof of Theorem~\ref{1.pl.sec.thm}  in positive characteristic.
Unfortunately,  the irreducibility of ${\mathcal H}_{d,g;d}$ is open even  over $\c$ and Hurwitz schemes are known to be more complicated in  characteristic $p>0$; cf.\ \cite{MR1837111}.
\end{defn}
 
Our aim is to show that an open dense subset of 
${\mathcal H}_{d,g;d}$ consists of maps whose ramification is as simple as possible. We use the discriminant as the relevant  measure.

\begin{prop}\label{hurw.open.prop}
 Let $k$ be a perfect field of characteristic $p\neq 2$. 
There is an open, dense subset ${\mathcal H}^\circ_{d,g;d}\subset {\mathcal H}_{d,g;d}$ such that for every $(\bar\pi:\bar B\to \p^1)\in {\mathcal H}^\circ_{d,g;d}$
\begin{enumerate}
\item $\bar\pi$ is separable,
\item $\mfrd_{\infty}(\bar\pi)=d-1$ if $p\nmid d$,  
\item $\mfrd_{\infty}(\bar\pi)=d$ if $p\mid d$ and
\item $\mfrd_{b}(\bar\pi)\leq 1$ for every $b\in B=\bar B\setminus\{b_{\infty}\}$. 
\end{enumerate}
\end{prop}

Proof. The properties (1--4) are all open, hence it remains to show that every
$\bar\pi:\bar B\to \c\p^1$ has a small deformation with these properties.

First we deal with (1). Let $p^e$ be the degree of inseparability of $\pi$. 
We can then factor $\bar\pi$ as
$$
\bar\pi: B\stackrel{\bar\pi^s}{\longrightarrow} \p^1\stackrel{F^e}{\longrightarrow} \p^1
$$
where $F^e$ is the $e$th power Frobenius given by
$(u{:}v)\mapsto  \bigl(u^{p^e}: v^{p^e}\bigr)$. The latter has separable deformations, for example  $(u{:}v)\mapsto  \bigl(u^{p^e}+tuv^{p^e-1}: v^{p^e}\bigr)$. Composing it with $\bar\pi^s$ gives a separable deformation
of $\bar\pi$.

It remains to prove that if a separable morphism $\bar\pi_0:B_0\to \p^1$
does not satisfy the conditions (2--4) then it has a  deformation
$\bar\pi_t:B_t\to \p^1$ with smaller discriminant.
A direct application of \cite[Thm.4.1]{MR0352540} shows that the latter is a local question at the ramification points; the relevant definitions and results are recalled in Paragraph~\ref{glob.curve.maps.say}.
Thus it remains to discuss how to lower the discriminant by local deformations.

We start with the ramification point at infinity.
 Choose local coordinates  $v$ at $b_\infty\in \bar B$ and
$u$ at $\infty\in \p^1$. Then $\bar\pi_0$ is given by a power series
$$
\phi_0(v)= a_d v^d+a_{d+1} v^{d+1}+(\mbox{higher terms}),
$$
where $a_d\neq 0$. Choose a deformation of it over $\spec k[[t]]$ given by
$$
\Phi(v,t)= a_d v^d+(t+a_{d+1}) v^{d+1}+(\mbox{higher terms}).
\eqno{(\ref{hurw.open.prop}.5)}
$$
The ramification index is $d$ for every $t$, thus we stay in ${\mathcal H}_{d,g;d}$. 
If $t\neq -a_{d+1}$ then either  $da_d$ or  $(d+1)(t+a_{d+1})$ is nonzero,
thus  $\mfrd_{\infty}(\pi_t)\leq d$ and (2--3) hold.

Choosing local coordinates at a ramification point in $B$, 
$\bar\pi_0$ is given by a power series
 $\phi_0(v)=a_nv^n+....+a_mv^m+...$  where $\mult \phi'_0(v)=m-1$ (thus $p\nmid m$). We choose the deformation 
$\Phi(v, t)=tv^2+\phi_0(v)$. Then
$$
\tfrac{\partial \Phi}{\partial v}=v\bigl(t+ma_mv^{m-2}+(\mbox{higher terms})\bigr).
$$
Here $t+ma_mv^{m-2}+(\mbox{higher terms})$ vanishes to order $m-2$ for $t=0$, so to order $\leq m-2$ nearby.   (In fact we have only simple ramification if $p\nmid m-2$.)
The $(v=0)$ branch is smooth and meets the previous branch  only at $t=0$. 
So we lowered the coefficients in $\mfrd$ from $m-1$ to $\leq m-2$. \qed
\medskip

{\it Note \ref{hurw.open.prop}.6.} The second part of the above argument 
works to simplify the ramification of any separable map
$B\to C$ between smooth, projective curves. By contrast, if
$g(C)\geq 2$ then a purely inseparable map $C\to C$ does not have separable deformations. This follows from the Hurwitz formula.

Comparing the bounds (\ref{pi1.onto.crit.lem}.1) with (\ref{hurw.open.prop}.2--3) gives the following.

\begin{lem} \label{pi1.onto.crit.lem.c} Let $k$ be a    field of characteristic
 $p\neq 2$. Then every  morphism   $(\pi:B\to \a^1)$  in 
${\mathcal H}^\circ_{d,g;d}$ is $\pi_1$-surjective.  \qed
\end{lem}

{\it Example \ref{pi1.onto.crit.lem.c}.1.} If $\chr k=2$ then
every degree 2, separable  morphism 
$\a^1\to \a^1$ is  \'etale. So none of the maps in
${\mathcal H}^\circ_{2,0;2}$ are $\pi_1$-surjective.

\begin{say}[Globalizing local deformations] \cite{MR0352540}
\label{glob.curve.maps.say}
Let $g:C\to B$ be a separable morphism of smooth projective curves with ramification points $c_i\in C$. Set $b_i=g(c_i)$.
Informally, we claim that  deformations of the local morphisms
$g_i:(c_i, C)\to (b_i, B)$ can be globalized. 

To make this assertion precise,  let 
$\hat g_i:\hat C_i\to \hat B_i$ denote  the completion of $g$ at $c_i$.
After choosing local coordinates $u_i$ at $b_i$ and $v_i$ at $c_i$
$\hat g_i$ is equivalent to an injective ring map
$k[[u_i]]\to k[[v_i]]$ given by $u_i\mapsto \phi_i(v_i)$.

For every $i$ let
 $\hat G_i:\hat {\mathbf C}_i\to \hat {\mathbf B}_i$
be a flat deformation of $\hat g_i:\hat C_i\to \hat B_i$. 
Equivalently,  a ring map
$k[[u_i,t]]\to k[[v_i,t]]$ given by $u_i\mapsto \Phi_i(v_i,t)$ where
$\Phi_i(v_i,0)=\phi_i(v_i)$. The main result is the following.

\medskip

{\it Theorem \ref{glob.curve.maps.say}.1.}  
 There is a pointed curve $(0,D)$ and a smooth, projective morphism
$G: {\mathbf C}\to B\times D$ such that
\begin{enumerate}
\item[(a)]  $\bigl(G_0: {\mathbf C}_0\to \{0\}\bigr)\cong \bigl(g:C\to B\bigr)$ and
\item[(b)]  the completion of $G$ at $(c_i, 0)$ is isomorphic to 
$\hat G_i:\hat {\mathbf C}_i\to \hat {\mathbf B}_i$
 for every $i$.
\end{enumerate}

Sketch of proof. Over $\c$ an argument goes back to Riemann.
We first construct $G: {\mathbf C}\to B\times D$ as a topological branched cover and then use Riemann's existence theorem to show that 
${\mathbf C} $ can be endowed with a unique complex structure such that
$G$ becomes holomorphic.

This argument is harder to do in full generality, and a theory
of deformations of morphisms was worked out in \cite{MR0352540}.
As written, it treats morphisms $g:X\to Y$ of complex manifolds in arbitrary dimension, but the arguments work in all characteristic without changes.
For us the relevant result is \cite[Thm.4.1]{MR0352540}, which says that
every deformation of the formal neighbourhood of the ramification locus of $g$ extends to a deformation   of $(g:X\to Y)$ if
\begin{enumerate}\setcounter{enumi}{1}
\item[(c)] $H^1(X, T_X)\to H^1(X, g^*T_Y)$ is surjective and
\item[(d)] $H^2(X, T_X)\to H^2(X, g^*T_Y)$ is injective.
\end{enumerate}
In the case of curves the $H^2$ are automatically zero and, by Serre duality, (c) is equivalent to the injectivity of
$$
H^0(X, \omega_X\otimes g^*\omega_Y)\to H^0(X, \omega_X^2).
$$
The latter holds if $g^*\omega_Y\to  \omega_X$ is nonzero, that is, when $g$ is separable. \qed
\end{say}

\section{Fundamental index in flat families}\label{fund.ind.def.sec}

In this section we work over a field $k$ of characteristic $p>0$.
All statements hold in  characteristic 0, but they only give a more convoluted proof of   Theorem~\ref{H1.surj.fund.sols.thm}.

\begin{prop}\label{fund.pulb.thm}  Let 
$X$ be a connected $k$-scheme, $\bar C$  a geometrically connected and geometrically reduced curve over $k$ and $D$ a  divisor on $\bar C^{\rm sm}$, the smooth locus of $\bar C$.
Let $\pi:\bar {\mathcal G}\to X\times \bar C$ be a finite, flat  morphism.
For $x\in X$ by base change we get $\pi_x:\bar G_x\to \bar C_x$. 
Assume the following.
\begin{enumerate}
\item The fibers of the projection $\bar {\mathcal G}\to X$ are geometrically connected and geometrically reduced curves.
\item $\bar G_x$ is smooth along $\pi_x^{-1}D_x$ for every $x\in X$.
\item $\pi_x$ has ramification index $e$ at every point of $\pi_x^{-1}D_x$  for every $x\in X$.
\end{enumerate}
Then   the prime-to-$p$ part of
 $\find\bigl(\pi_x^*D_x\bigr)$ is  independent of $x\in X$.
\end{prop}

Proof.  Write $D=mD'$ where $D'$ is primitive. Then $\tfrac1{me}\pi_x^*D_x$
is primitive for every $x$. 
By (\ref{prim.fund.defn}.1)
$$
\find\bigl(\pi_x^*D_x\bigr)\cdot 
\ord\bigl(\tfrac1{me}\pi_x^*D_x\bigr)=me.
\eqno{(\ref{fund.pulb.thm}.4)}
$$
We check in Proposition~\ref{tors.ord.loc.const} that the prime-to-$p$ part of
 $\ord\bigl(\tfrac1{me}\pi_x^*D_x\bigr)$ is independent of $x\in X$.  
Then (\ref{fund.pulb.thm}.4) shows that  the prime-to-$p$ part of
 $\find\bigl(\pi_x^*D_x\bigr)$ is also independent of $x\in X$. 
\qed

\begin{prop} \label{tors.ord.loc.const}
Let $X$ be an irreducible scheme with generic point $x_g\in X$.
Let $\tau:  \bar{\mathcal C}\to X$ be a 
flat, projective morphism whose fibers are geometrically connected and geometrically reduced curves. Let $L$ be a line bundle on $\bar{\mathcal C} $
and assume that $\ord(L_x)$ is finite for every $x\in X$. 
Then $\ord\bigl(L_{x_g}\bigr)=\ord\bigl(L_{x}\bigr)\cdot p^{c(x)}$  for some $c(x)\geq 0$.
\end{prop}

Proof. Write $\ord\bigl(L_{x_g}\bigr)=np^c$ 
where $p\nmid n$ and replace $L$ by $L^{p^c}$. 
We can thus assume that $\chr(k(x))\nmid n:=\ord\bigl(L_{x_g}\bigr)$.

As in Paragraph~\ref{jac.say.2}, we have a universal family of Jacobians
$\jac\bigl(\bar {\mathcal C}/X\bigr)\to X$ and $L$ gives
 a section 
$$
\sigma_L: X\to\jac\bigl(\bar {\mathcal C}/X\bigr).
$$
For any $d$ let  ${\mathcal T}[d]\subset \jac\bigl(\bar {\mathcal C}/X\bigr)$ denote the $d$-torsion subgroup.  Let $\sigma_0:X\to \jac\bigl(\bar {\mathcal C}/X\bigr)$ denote the zero-section.
Note that ${\mathcal T}[d]\to X$ 
 is \'etale
over $x\in X$ whenever $\chr(k(x))\nmid d$
(see, for example \cite[p.64]{mumf-abvar}).

By our choice  the image of 
$\sigma_L $ lies in 
${\mathcal T}[n]$.   Set $m:=\ord\bigl(L_{x}\bigr)$ and consider the sections
$\sigma_{L^m}$ and $ \sigma_0$. By assumption
$ \sigma_{L^m}(x)=\sigma_0(x)$. Since  ${\mathcal T}[n]\to X$ 
 is \'etale, this implies that $ \sigma_{L^m}=\sigma_0$. Thus
$\ord\bigl(L_{x_g}\bigr)=m=\ord\bigl(L_{x}\bigr)$. \qed

\begin{say}[Relative Jacobians] \label{jac.say.2}
Let $\tau:  \bar{\mathcal C}\to X$ be a flat, proper morphism 
whose fibers are geometrically reduced and geometrically connected curves.
(In particular,  $H^0(\bar C_x, \o_{\bar C_x})\cong k(x)$ and 
the dimension of $H^1(\bar C_x, \o_{\bar C_x})$, which is also the dimension of
$\jac(\bar C_x)$, is locally constant on $X$.) 
Then the Jacobians of the fibers form a flat family
$$
 \jac\bigl(\bar{\mathcal C}/X\bigr)\to X.
$$
This follows from the theory of Picard varieties outlined in
\cite{FGA} and a more detailed treatment of this case  can be found  in 
\cite[Chap.9]{blr}.

\end{say}

\section{Endomorphisms of Pell surfaces}\label{endom.sec}

Let us start by writing down some
 endomorphisms of $S_g$.

\begin{say}[Examples of endomorphisms] Let $\pi:S_g\to B$ be a Pell surface.
 The automorphisms of
$\bigl(B, (g=0)\bigr)$  form a  group that we denote by $\aut(B,g)$. We call these the {\it base automorphisms.} This group is infinite only in cases (\ref{special.say}.1--2) and (\ref{special.say}.4).
Otherwise, in characteristic 0 this group is cyclic, but in
positive  characteristic we can have the larger group of all affine linear transformations of $\f_q$. 

Let  $\sigma\in \aut(B, g)$.
Then $\sigma(g)=cg$ for some constant $c$. If $c$ is a square in the base field then we can lift $\sigma$ to an 
automorphism of $S_g$ (in 2 ways) by setting
$\sigma':(x,y,u)\mapsto  \bigl(x,c^{-1/2}y,\sigma(u)\bigr)$.

If $\Sigma$ is a section of $\pi:S_g\to B$ then {\it translation by $\Sigma$} is an automorphism of
$S_g$. Explicitly, if $\Sigma=\bigl(s_x(u), s_y(u), u\bigr)$ is a section  then we get
$$
(x,y,u)\mapsto  \bigl(s_x(u)x+gs_y(u)y, s_x(u)y+s_y(u)x,u\bigr)
$$
The translation subgroup of $\aut(S_g)$ is  either  $\z/2$ or $\z+\z/2$.

For any $n\in \z$ the {\it $n$th power map} is an endomorphism of $S_g$.
For $n=-1$  we get an automorphism, the inverse map, which is
$(x,y,u)\mapsto   (x,-y,u)$. For $n\geq 1$ the map is given as
$(x,y,u)\mapsto \bigl(x_n, y_n, u\bigr)$ where
$$
\begin{array}{lcl}
x_n& = & \sum_{i=0}^{\rdown{n/2}} \binom{n}{2i}x^{n-2i}y^{2i} g(u)^i\qtq{and}\\
y_n& = & \sum_{i=0}^{\rdown{n/2}} \binom{n}{2i+1}x^{n-2i-1}y^{2i+1}g(u)^i.
\end{array}
$$
The translations and the inverse map generate a subgroup of $\aut(S_g,\pi)\subset \aut(S_g)$ which is  either   $\z/2+\z/2$ or $D_{\infty}+\z/2$,
where $D_{\infty}$ denotes the infinite dihedral group.
\end{say}

\begin{thm}\label{endom.thm} Let $\pi:S_g\to B$ be a   Pell surface  other than the exceptional ones 
(\ref{special.say}.1--4).
\begin{enumerate}
\item  There is an exacxt sequence $1\to \aut(S_g,\pi)\to \aut(S_g)\to \aut(B,g)$.
\item Every dominant, separable endomorphism  $\phi$ of $S_g$ can be 
written (non-uniquely) as
$\phi=\phi^{\rm end}\circ \phi^{\rm aut}$ where $\phi^{\rm aut}$ is an automorphism and
$\phi^{\rm end}$ is the $n$th power map for some $n\geq 1$.
\end{enumerate}
\end{thm}

Proof. Since $\phi$ is separable, pulling back by $\phi$ 
gives an injection (hence an isomorphism) on $m$-canonical forms with log poles. 
If the log Kodaira dimension of $S_g$ is 1, then   $\pi:S_g\to B$ is the Iitaka fibration, 
hence $\phi$ commutes with $\pi$, up to an element of $\aut(B,g)$.
After precomposing with the inverse of the latter, we may thus assume  that
$\phi$ commutes with $\pi$, hence it induces an endomorphism of the generic
fiber. The geometric generic fiber is isomorphic to $\gm$, hence its
endomorphisms  are the translations composed with   power maps.
If this extends to an endomorphism of $S_g$ then the translation must map the identity section  to another  section. Hence only translations by a section extend to automorphisms of $S_g$.

By  (\ref{pell.surf.lem.1}.5) this leaves open only the case
$S_2=\bigl(x^2-(u^2-1)y^2=1\bigr)$. 
The above arguments apply to those endomorphisms that commute with $\pi$, thus it remains to show that there are no other
dominant, separable endomorphisms, up to the action of $\aut(\a^1, g)$.

\begin{say}[Intersection points of sections on $S_g$]\label{int.of.sec.S2}
Consider a Pell surface $S_g$ over an algebraically closed field. 
 Let $I_g\subset S_g$ denote the set of
all intersection points of distinct sections.

\medskip
{\it Claim \ref{int.of.sec.S2}.1.}  Let $S_g$ be a Pell surface with a nontrivial section. Then
\begin{enumerate}
\item[(a)] $I_g$ is Zariski dense in $S_g$ and
\item[(b)]  every point of $I_g$ is contained in
infinitely many  sections. 
\end{enumerate}
\medskip

Proof. 
Let $F_b\subset S_g$ be an irreducible fiber. Fixing a value of  $\sqrt{g(b)}$,
the map $T_b$ defined in (\ref{pell.surf.defn}.3) gives a  group homomorphism
on the group of sections 
$$
\tau_b:\Sigma\mapsto \Sigma\cap F_b\mapsto T_b(\Sigma\cap F_b)\in \gm.
$$
Assume now that there are infinitely many sections.
Let us call $F_b$ a {\it cyclotomic fiber} if the image of $\tau_b$ is finite.
If $\Sigma_1=(x_1, y_1)$ is a fundamental section then
$F_b$ is cyclotomic iff $ x_1(b)+y_1(b)\sqrt{g(b)}$ is a root of unity.
If the root of unity has order $r$ then  the different sections meet
$F_b$ in $r$ points if $r$ is even and  $2r$ points if $r$ is odd. 

Since $x_1+y_1 \sqrt{g}$ is non-constant, there are infinitely many cyclotomic fibers and the order of   $ x_1(b)+y_1(b)\sqrt{g(b)}$ is unbounded. \qed
\medskip

We can be even more precise for  $g=u^2-1$. 
The fundamental section is then $u\mapsto (u,1)$ thus  $F_b$ is a  cyclotomic fiber iff $b+\sqrt{b^2-1}=\zeta$ is a root of unity. Thus $b=\frac12(\zeta+\zeta^{-1})$ and  we obtain the following.
\medskip

{\it Claim \ref{int.of.sec.S2}.2.} If $\chr k=0$ then the  cyclotomic fibers on $S_2$ are exactly the ones lying over the points of $R_{\infty}:=\{\cos(2\pi \alpha): \alpha\in \q\}$. If $\chr k=p>0$ then the  cyclotomic fibers on $S_2$ are exactly the ones lying over the points of  $\bar\f_p$. 
\qed
\medskip

Since a dominant morphism between Pell surfaces can map only finitely many sections to the same section, we get the following.

\medskip
{\it Claim \ref{int.of.sec.S2}.3.} Let $\phi:S_g\to S_h$ be a dominant morphism between Pell surfaces. Assume that $S_g$ has a nontrivial section.
Then  $\phi(I_g)\subset I_h$. \qed
\end{say}

\begin{say}[Endomorphisms of $S_2$ in characteristic $0$] \label{end.S2.char.0.say} 
As we computed in (\ref{pell.surf.lem.1}.7), $x^{-1}dy\wedge du$ is the unique (up to scalar)
2-form with log poles at infinity on $S_2$.
Let $\phi:S_2\to S_2$ be a dominant endormorphism. Then $\phi^*\bigl(x^{-1}dy\wedge du\bigr)$ is also a 2-form with log poles at infinity, hence a scalar multiple of $x^{-1}dy\wedge du$. Thus $\phi$ is \'etale.  We also know that $\phi$  maps affine lines to (possibly singular) affine lines and
an \'etale morphism  $\a^1\to \a^1$ is an isomorphism. 
Thus if $\Sigma$ is  a section, then (with finitely many possible exceptions due to vertical lines) $\Sigma':=\phi(\Sigma)$ is another section  and 
$\phi$ gives an isomorphism $\Sigma\to \Sigma'$. 
 Thus, in the diagram below,
3 of the maps are isomorphisms, hence so is the bottom arrow which we call $\tau_{\Sigma}$.
$$
\begin{array}{ccc}
\Sigma  & \stackrel{\phi}{\to} & \Sigma'\\
\pi\downarrow\hphantom{\pi} && \hphantom{\pi}\downarrow \pi\\
\a^1  &\stackrel{\tau_{\Sigma}}{\to} &  \a^1.
\end{array}
\eqno{(\ref{end.S2.char.0.say}.1)}
$$
Thus $\tau_{\Sigma}(u)=a_{\Sigma}u+b_{\Sigma}$ for some $a_{\Sigma}, b_{\Sigma}$. Note further that
(\ref{int.of.sec.S2}.3) implies that 
$\tau_{\Sigma}$  maps  $R_{\infty}$ to itself. So
 $\tau_{\Sigma}(u)=\pm u$ by  (\ref{end.S2.char.0.say}.2).

After precomposing with $u\to -u$ if necessary, 
we may thus assume that 
 there are infinitely many sections  $\{\Sigma_n: n\in I\}$  on which 
$\phi$ commutes with $\pi$.
Let now $F_c$ be a non-cyclotomic fiber. Then 
$\{\Sigma_n\cap F_c: n\in I\}$ is an infinite  subset of $F_c$
that is mapped to $F_c$ by $\phi$.  Thus $F_c\cap \phi(F_c)$ is infinite, hence $\phi(F_c)=F_c$. This shows that  $\phi$
commutes with the projection $\pi$.  This completes the proof of
Theorem~\ref{endom.thm} for $S_2$ in characteristic $0$.\qed
\medskip

{\it Lemma \ref{end.S2.char.0.say}.2.}  Assume that $p(x)=ax+b$ maps  $R_{\infty}$ to itself. If $\chr k=0$ then
$b=0$ and $a=\pm 1$.

Proof.
By assumption there is an $n$ such that 
$$
p(1), p(-1)\in R_n:=\bigl\{\cos(2\pi \alpha): \alpha\in \tfrac1{n}\z\bigr\}.
$$ Thus
$a, b\in \q(\zeta_n+\bar\zeta_n)$ for some $n$, and so  $p$ maps $R_n$ to itself
injectively. Since $R_n$ is finite, $p:R_n\to R_n$ is a bijection, hence
$p\bigl([-1,1]\bigr)= [-1,1]$. \qed
\end{say}

\begin{say}[Endomorphisms of $S_2$ in characteristic $p>0$] \label{end.S2.char.p.say} 
Let $\phi:S_2\to S_2$ be a dominant endomorphism of degree $d$. Then  $\phi$  maps affine lines to (possibly singular) affine lines.  
Thus $\phi$ maps $I_g$ to itself, so $\phi$  is defined over $\bar\f_p$ by  (\ref{end.S2.char.p.say}.3), hence over 
  a finite field $\f_q$ for some $q=p^c$.
As before we get a
commutative diagram
$$
\begin{array}{ccc}
\Sigma  & \stackrel{\phi}{\to} & \Sigma'\\
\pi\downarrow\hphantom{\pi} && \hphantom{\pi}\downarrow \pi\\
\a^1_t  &\stackrel{\tau_{\Sigma}}{\to} &  \a^1_u,
\end{array}
\eqno{(\ref{end.S2.char.p.say}.1)}
$$
where we only know that $\deg \tau_{\Sigma}\leq d$.  So
 $\tau_{\Sigma}$ is a degree $\leq d$ polynomial over $\f_q$.  Since the latter form a finite set, there are infinitely many sections  $\{\Sigma_n: n\in I\}$   with the  
same $\tau_{\Sigma}$; call this common map $\tau$. 
Let now $F_c$ be a non-cyclotomic fiber. (This always exists after a transcendental base field extension.) Then 
$\{\Sigma_n\cap F_c: n\in I\}$ is an infinite  subset of $F_c$
that is mapped to $F_{\tau(c)}$ by $\phi$.  Thus $F_{\tau(c)}\cap \phi(F_c)$ is infinite, hence $\phi(F_{\tau(c)})=F_c$. This shows that  $\phi$
sits in a commutative diagram
$$
\begin{array}{ccc}
S_2  & \stackrel{\phi}{\to} & S_2\\
\pi\downarrow\hphantom{\pi} && \hphantom{\pi}\downarrow \pi\\
\a^1_t  &\stackrel{\tau}{\to} &  \a^1_u.
\end{array}
\eqno{(\ref{end.S2.char.p.say}.2)}
$$
Thus $\phi$ factors through the Pell surface
$$
S_h:=\bigl(x^2-\bigl(\tau(t)^2-1\bigr)y^2=1\bigr)
$$
This is only possible if  the curve $C_g$ has genus 0 for
$g=\tau(t)^2-1$. That is, when $g$ has exactly 2 roots of odd multiplicity. 
As we discussed in Example~\ref{all.mult.roots}, this only happens when
$\tau(t)=T_n(t)$ for some $n$, where  $T_n(t)$ is the $n$th Chebyshev polynomial of the first kind as in  (\ref{u2-1.exmp}.4).

These examples give the very interesting rational maps
$$
\phi_n:S_2\map S_2\qtq{given by}
(x, y, t)\mapsto \Bigl(x, U_n(t)^{-1}\cdot y, T_n(t)\Bigr).
$$
 These maps are, however,  not defined along $U_n(t)=0$.
This completes the proof of
Theorem~\ref{endom.thm} for $S_2$ in characteristic $>0$.\qed\medskip
\medskip

{\it Lemma \ref{end.S2.char.p.say}.3.}  Let $X, Y$ be  $K$-varieties,
$L/K$ a field extension and  $\phi:X_L\to Y_L$ a morphism. 
Assume that there is a Zariski dense set $S\subset X(K)$ such that
$\Phi(S)\subset Y(K)$. Then $\phi$ is defined over $K$. \qed
\medskip

\end{say}

\begin{ques} It is natural to ask if all morphisms $S_h\to S_g$ between Pell surfaces
are  compositions of endomorphisms and of the base change maps
$S_{g\circ q}\to S_g$. Our methods settle this if the log Kodaira dimension of $S_g$ is 1, or if  $S_h$ contains infinitely many affine lines. 
The remaining step is to understand  all morphisms $S_h\to S_2$
when $S_h$ contains only the obvious affine lines.  
\end{ques}


\def\cprime{$'$} \def\cprime{$'$} \def\cprime{$'$} \def\cprime{$'$}
  \def\cprime{$'$} \def\dbar{\leavevmode\hbox to 0pt{\hskip.2ex
  \accent"16\hss}d} \def\cprime{$'$} \def\cprime{$'$}
  \def\polhk#1{\setbox0=\hbox{#1}{\ooalign{\hidewidth
  \lower1.5ex\hbox{`}\hidewidth\crcr\unhbox0}}} \def\cprime{$'$}
  \def\cprime{$'$} \def\cprime{$'$} \def\cprime{$'$}
  \def\polhk#1{\setbox0=\hbox{#1}{\ooalign{\hidewidth
  \lower1.5ex\hbox{`}\hidewidth\crcr\unhbox0}}} \def\cdprime{$''$}
  \def\cprime{$'$} \def\cprime{$'$} \def\cprime{$'$} \def\cprime{$'$}
\providecommand{\bysame}{\leavevmode\hbox to3em{\hrulefill}\thinspace}
\providecommand{\MR}{\relax\ifhmode\unskip\space\fi MR }
\providecommand{\MRhref}[2]{%
  \href{http://www.ams.org/mathscinet-getitem?mr=#1}{#2}
}
\providecommand{\href}[2]{#2}

\bigskip

\noindent  Princeton University, Princeton NJ 08544-1000

\email{kollar@math.princeton.edu}

\end{document}